\documentclass[preprint,12pt]{elsarticle}

\usepackage{amsmath, amsthm, amscd, amsfonts, amssymb}
\def\red#1{{\textcolor{Red}{#1}}}
\usepackage[usenames,dvipsnames]{color}

\newtheorem{theorem}{Theorem}[section]

\newtheorem{lemma}[theorem]{Lemma}
\theoremstyle{definition}
\newtheorem{definition}[theorem]{Definition}
\newtheorem{remark}[theorem]{Remark}

\def\st#1#2{{\mathrel{\mathop{#2}\limits_{#1}}{}\!}}
\def\vsp{\vspace*{2mm}\\ }

\journal{Journal Functional Analysis}

\begin{document}

\begin{frontmatter}

\title{Solutions for nonlinear Fokker--Planck equations with~measures as initial data and~McKean-Vlasov equations}

\author[VBarbu]{Viorel Barbu}
\author[MRockner]{Michael R\"ockner}

\address[VBarbu]{Octav Mayer Institute of Mathematics of the Romanian Academy,     Ia\c si, Romania.  Email: vbarbu41@gmail.com}
\address[MRockner]{Fakult\"at f\"ur Mathematik, Universit\"at Bielefeld,  D-33501 Bielefeld, Germany,\\   and Academy of Mathematics 
 and System Sciences, CAS, Beijing, China\\  Email: roeckner@math.uni-bielefeld.de}

\begin{abstract}
One proves the existence and uniqueness of a generalized (mild) solution for the nonlinear Fokker--Planck equation (FPE)
$$\begin{array}{ll}
u_t-\Delta (\beta(u))+{\rm div}(D(x)b(u)u)=0,&t\ge0,\ x\in\mathbb{R}^d,\ d\ne2,\vsp
u(0,\cdot)=u_0,\mbox{in }\mathbb{R}^d,\end{array}$$
where $u_0\in L^1(\mathbb{R}^d)$, $\beta\in C^2(\mathbb{R})$ is a  nondecreasing  function,
$b\in C^1$, bounded, $b\ge0$, $D\in {L^\infty}(\mathbb{R}^d;\mathbb{R}^d)$, ${\rm div}\,D\in L^2(\mathbb{R}^d)+L^\infty(\mathbb{R}^d),$ with ${({\rm div}\, D)^-}\in L^\infty(\mathbb{R}^d)$,  $\beta$ strictly increasing,  if $b$ is not constant.   Moreover, $t\to u(t,u_0)$ is a semigroup of contractions in $L^1(\mathbb{R}^d)$, which leaves invariant the set of  probability density functions in $\mathbb{R}^d$. If ${\rm div}\,D\ge0$, $\beta'(r)\ge a|r|^{\alpha-1}$,    and $|\beta(r)|\le C r^\alpha$, $\alpha\ge1,$  $d\ge3$, then $|u(t)|_{L^\infty}\le Ct^{-\frac d{d+(\alpha-1)d}}\ |u_0|^{\frac2{2+(m-1)d}},$ $t>0$, and,  if $D\in L^2(\mathbb{R}^d;\mathbb{R}^d)$,  the existence extends to initial data $u_0$ in the space $\mathcal{M}_b$  of bounded  measures in $\mathbb{R}^d$.  As a consequence for arbitrary initial laws, we obtain weak  solutions to a class of  McKean-Vlasov SDEs with coefficients which have singular dependence on  the time marginal laws.
\end{abstract}


\end{frontmatter}

\section{Introduction}\label{s1}

Here, we consider the nonlinear  parabolic  equation (FPE)
\begin{equation}\label{e1.1}
\begin{array}{ll}
u_t-\Delta(\beta(u))+{\rm div}(Db(u)u)=0,\ \forall (t,x)\in[0,{\infty})\times{\mathbb{R}}^d,\\
 u(0,\cdot)=\mu,\ \mbox{in }{\mathbb{R}}^d,\end{array}\end{equation}
where $\mu$ is a bounded measure on ${\mathbb{R}}^d$ and the functions $\beta:{\mathbb{R}}\to{\mathbb{R}}$,\break  $D:{\mathbb{R}}^d\to{\mathbb{R}}^d$, $b:{\mathbb{R}}\to{\mathbb{R}}$, are assumed to satisfy the following hypotheses
\begin{itemize}
	\item[(i)] {\it $\beta:{\mathbb{R}}\to{\mathbb{R}}$ is  a  monotonically nondecreasing $C^1$-function, $\beta(0)=0$.}
	\item[(ii)] $D\in  L^{\infty}({\mathbb{R}}^d;{\mathbb{R}}^d),$ ${\rm div}\,D\in L^1_{\rm loc}(\mathbb{R}^d),$ $({\rm div}\ D)^-\in L^{\infty}({\mathbb{R}}^d).$
	\item[(iii)] {\it$b\in C^1({\mathbb{R}}),$ $b$ bounded, nonnegative, and $b\equiv{\rm const.}$, if $\beta$ is not  strictly increasing.}
\end{itemize}
In statistical physics, the nonlinear Fokker--Planck equation (NFPE) \eqref{e1.1} models  anomalous diffusion processes.   As a matter of fact, \eqref{e1.1} is an extension of the classical Smoluchowski equation.

It should be also mentioned that, as in the classical linear theory, NFPE  \eqref{e1.1} is related to the McKean-Vlasov stochastic differential equation
\begin{equation}\label{e1.2}
\begin{array}{l}
dX(t)=D(X(t))b(u(t,X(t)))dt+\frac1{\sqrt{2}}\left(\frac{\beta(u(t,X(t))}{u(t,X(t))}\right)^{\frac12}dW(t),\ t\ge0,\\
X(0)=\xi_0,\end{array}\end{equation}where  $\frac{\beta(0)}0:=\beta'(0).$  
More precisely, if $u$ is a Schwartz distributional solution to \eqref{e1.1} such that $u:[0,{\infty})\to L^1({\mathbb{R}}^d)$ is $\sigma(L^1,C_b)$-continuous with \mbox{$u(0)=u_0\,dx$,} $u_0\in L^1\cap L^{\infty}$, $u\ge0,$ and
$$\int_{{\mathbb{R}}^d} u(t,x)dx=1,\ \forall  t\ge0,$$then there exists a (probabilistically) weak solution process $X$ to SDE \eqref{e1.2} such that $u$ is the probability density of its time marginal laws.   (See \mbox{\cite{4}--\cite{5}}.)  Here we would like to note that there is a vast literature on McKean--Vlasov SDEs and on Fokker--Planck equations (see \cite{15'} and  in particular \cite{11'}, Section 6.7 (iii), respectively, and the references therein). But the dependence of the coefficients on the solution in these references is much more regular than in our so called "Nemytskii case" (see \cite{5}, Section 2, for details).

The existence theory in ${\mathbb{R}}^d$ for   Fokker--Planck equations of the form \eqref{e1.1} in ${\mathbb{R}}^d$, $d\ge1$, and weak existence for the corresponding McKean--Vlasov equations \eqref{e1.2}  was studied in \cite{4}-\cite{5}, \cite{6a}, via  nonlinear semigroup theory in  $L^1({\mathbb{R}}^d)$.  Here, we shall extend this approach to the case where the initial data is a  bounded   measure on ${\mathbb{R}}^d$.   In this case, our existence result extends those in earlier work of  H.~Brezis and A.~Friedman \cite{5a} and by  M. Pierre \cite{13a}, who studied the case $D\equiv0$. One main ingredient of the proof is the $L^1-L^{\infty}$ regularity of the solution of \eqref{e1.1}, extending  previous results by Ph.~Benilan \cite{8a}, A.~Pazy \cite{14},   L.~Veron \cite{15}.

A further contribution of this paper is its  applications to McKean-Vlasov SDEs (see Section~\ref{s4}) which in turn was one motivation of this work.

The structure of the paper is as follows.
Sections \ref{s2} and \ref{s3} contain  our main exis\-tence result for \eqref{e1.1} for $\mu=u_0\,dx$, $u_0\in L^1$.
In Section \ref{s5}, we prove the $L^1-L^{\infty}$ regularizing effects of the generalized solutions to \eqref{e1.1} on the initial data (see Theorem \ref{t5.1}).
Section \ref{s6} is devoted to the existence for the Radon measure initial data and Section \ref{s4} to applications to McKean-Vlasov SDEs.

\medskip \noindent{\bf Notation.} For each $1\le p\le{\infty}$, $L^p({\mathbb{R}}^d)$, briefly denoted   $L^p$, is the space of Lebesgue integrable functions $u$ on ${\mathbb{R}}^d$   with the standard norm $|\cdot|_p$.   We~denote by $W^{k,p}({\mathbb{R}}^d)$  the Sobolev space of all functions in $L^p$ with partial derivatives $D^\ell_j$ up to order $k$ in $L^p$,  where ${D^\ell_j}=\frac{\partial^\ell}{\partial x^\ell_j}$ is taken in the sense of Schwartz distributions ${\mathcal{D}}'({\mathbb{R}}^d)$.    Denote by $L^p_{\rm loc}$ the space of $L^p$-integrable functions on every compact set and $W^{k,p}_{\rm loc}=\{u\in L^p_{\rm loc}; \ D^\ell_ju \in L^p_{\rm loc},\ 1\le j\le d,\   \ell=1,2,...,k\}$, $1\le p\le{\infty}$.   We set $H^1=W^{1,2}({\mathbb{R}}^d)$ and $H^2=W^{2,2}({\mathbb{R}}^d)$. We denote by $\left<\cdot,\cdot\right>_2$ the scalar product of $L^2$ and by ${}_{H^{-1}}\left<\cdot,\cdot\right>_{H^1}$ the duality functional between $H^1$ and its dual space $(H^1)'=H^{-1}.$ 	
 $C_b({\mathbb{R}}^d)$ is the space of  all continuous and bounded functions \mbox{$u:{\mathbb{R}}^d\to{\mathbb{R}}$} endowed with the supremum norm $\|\cdot\|_{C_b}$, and by $C^1_b({\mathbb{R}}^d)$ the space of continuously differentiable bounded functions with bounded derivatives. By ${\mathcal{M}}_b({\mathbb{R}}^d)$ we denote the space of bounded Radon measures on ${\mathbb{R}}^d$. In the following, we shall simply denote   $C_b({\mathbb{R}}^d)$, $C^1_b({\mathbb{R}}^d)$ by $C_b$, $C^1_b$, respectively, and ${\mathcal{M}}_b({\mathbb{R}}^d)$ by ${\mathcal{M}}_b$. Denote by $\|\cdot\|_{{\mathcal{M}}_b}$ the norm of ${\mathcal{M}}_b$,
 $\|\mu\|_{{\mathcal{M}}_b}=\sup\{|\mu(\psi);\ |\psi|_{C_b}\le1\}.$
 We   set
	$${\mathcal{P}}_0({\mathbb{R}}^d)=\left\{\rho\in L^1,\ \rho\ge0,\ \int_{{\mathbb{R}}^d}\rho\,dx=1\right\}.$$	
	For $1<p<{\infty}$, we denote by $M^p({\mathbb{R}}^d)$ (simply written $M^p$) the Marcin\-kiewicz space of all (classes of) measurable functions $u:{\mathbb{R}}^d\to{\mathbb{R}}$ such that
	
	 $$\begin{array}{r}
	 \|u\|_{M^p}=\inf\Big\{{{\lambda}}>0;\displaystyle \int_K|u(x)|dx\le{{\lambda}}({\rm meas}\ K)^{\frac1{p'}}\\
	 \mbox{for all Borel sets $K\subset{\mathbb{R}}^d$}\Big\}<{\infty},\mbox{$\frac1p+\frac1{p'}=1.$}\end{array}$$
	
	 We recall that $M^p({\mathbb{R}}^d)\subset L^q_{\rm loc}({\mathbb{R}}^d)$   for \mbox{$1<q<p<{\infty}.$}
	
	Let $E(x)={\omega}_d|x|^{2-d}_d,$ $x\in{\mathbb{R}}^d$, $d\ge3$, be the fundamental solution of  the Laplace operator. (Here ${\omega}_d$ is the volume of the unit $d$-ball and $|\cdot|_d$ is the Euclidean norm of ${\mathbb{R}}^d$.) We recall that (see, e.g., \cite{6aa}) that, for $d\ge3,$
	\begin{equation}\label{e1.4}
E\in M^{\frac d{d-2}}({\mathbb{R}}^d),\ |\nabla E|_d\in M^{\frac d{d-1}}({\mathbb{R}}^d),
\end{equation}
and, for  $f\in L^1$, the solution $u$ to the equation
$-\Delta u=f\mbox{ in }{\mathcal{D}}'({\mathbb{R}}^d)$ is given by the convolution product $u=E*f$ and satisfies
\begin{eqnarray}
\label{e1.5}\|u\|_{M^{\frac d{d-2}}}&\le&\|E\|_{M^{\frac d{d-2}}}  |f|_1.
\\[1mm]
\label{e1.6}\|\nabla u\|_{M^{\frac d{d-1}}}&\le&\|\nabla E\|_{M^{\frac d{d-1}}}  |f|_1.
\end{eqnarray}(In this context, we cite also the work of H. Brezis and W. Strauss \cite{12a}.)

\section{Generalized solutions to NFPE \eqref{e1.1}}\label{s2}
\setcounter{equation}{0}
	
We shall treat first FPE \eqref{e1.1} for initial data $u_0\in L^1$.
	
	\begin{definition}\label{d2.1} \rm A continuous function $u:[0,{\infty}]\to  L^1$ is said to be a {\it generalized solution} (or {\it mild solution}) to equation \eqref{e1.1}  if, for each $T>0$,
\begin{equation}\label{e2.1}
		u(t)=\lim_{h\to0} u_h(t)\mbox{ in } L^1\mbox{ uniformly on each interval  $[0,T]$},\end{equation}
		where $u_h:[0,T]\to L^1$ is the step function, defined by the finite difference scheme	
 \begin{eqnarray}
 &&u_h(t)=u^i_h,\ \forall  t\in(ih,(i+1)h],\ i=0,1,...,N-1,\label{e2.2}\\
&&u^0_h=u_0,\ \beta(u^i_h)\in L^1_{\rm loc},\ u^i_h\in L^1,\ \forall  i=1,2,...,N,\ Nh=T,\label{e2.3}\\
&&u^{i+1}_h-h\Delta(\beta(u^{i+1}_h))+h\ {\rm div}(Db(u^{i+1}_h)u^{i+1}_h)=u^i_h,\ \mbox{in }{\mathcal{D}}'({\mathbb{R}}^d),\qquad
\label{e2.4}\end{eqnarray}
for all $i=0,1,....,N-1.$
\end{definition}
 Theorem \ref{t2.2} is our first main result.

 \begin{theorem}\label{t2.2} Let $d\ne2.$ Under Hypotheses {\rm(i), (ii), (iii)}, for each $u_0\in  \overline{D(A)}$, where $A$ and $ {D(A)}$ are defined in \eqref{e3.5a}, \eqref{e3.5aa} below and $\overline{D(A)}$ is the $L^1$-closure of $D(A)$, there is a unique generalized solution $u=u(t,u_0)$ to equation \eqref{e1.1}. Now assume, in addition, that
 	\begin{equation}
 	\label{e2.4prim}\beta\in C^2({\mathbb{R}}).
 	\end{equation} Then $\overline{D(A)}=L^1$ and, for every $u_0\in L^1\cap L^{\infty}$,
 	\begin{equation}\label{e2.4a}
 	|u(t)|_{\infty}\le
 	\exp(|({\rm div\,D})^-+|D||^{\frac12}_{\infty}t)|u_0|_{\infty},\ \forall  t>0.\end{equation}
 	If $u_0\in{\mathcal{P}}_0({\mathbb{R}}^d)$, then
 \begin{equation}
 u(t)\in {\mathcal{P}}_0({\mathbb{R}}^d),\ \forall  t\ge0,\label{e2.5} \end{equation}
 Moreover,  $t\to S(t)u_0=u(t,u_0)$ is a strongly continuous semigroup of nonlinear contractions from $L^1$ to  $L^1$, that is, $S(t+s)u_0=S(t)S(s)u_0$ for $0<s<t$, and
 \begin{equation}
 \label{e2.8}
 |S(t)u_0-S(t)\bar u_0|_1\le|u_0-\bar u_0|_1,\ \forall  t>0,\ u_0,\bar u_0\in L^1.
 \end{equation}
 If $u_0\in L^1\cap L^{\infty}$, then
  $u$ is a solution to \eqref{e1.1} in the sense of Schwartz distributions on $(0,{\infty})\times{\mathbb{R}}^d$, that is,
 \begin{equation}\label{e2.7}
\begin{array}{l}
\displaystyle  \int^{\infty}_0\int_{{\mathbb{R}}^d}(u({\varphi}_t+b(u)D\cdot\nabla{\varphi})+\beta(u)\Delta{\varphi})dt\,dx\vsp
\qquad+\displaystyle \int_{{\mathbb{R}}^d}u_0(x){\varphi}(0,x)dx=0,\ \forall {\varphi}\in C^{\infty}_0([0,{\infty})\times{\mathbb{R}}^d).\end{array}
 \end{equation}
\end{theorem}

  We also note   that, by \eqref{e2.4a}, equation \eqref{e2.7} is well defined for all ${\varphi}\in C^\infty_0([0,\infty)\times{\mathbb{R}}^d).$

\begin{remark}\label{r2.4}\rm
	It should be noted that the uniqueness of the solution $u$ given by Theo\-rem \ref{t2.2} is claimed in the class of generalized solutions and not in that of distributional solutions. The latter is true in some special cases   (see, e.g., \cite{6b},   \cite{7az}), but it is open in the general case we consider here.	
\end{remark}

\section{Proof of Theorem \ref{t2.2}}\label{s3}
\setcounter{equation}{0}

The idea of the proof is to associate with equation \eqref{e1.1} an $m$-accretive operator $A$ in $L^1$ and so to reduce \eqref{e1.1} to the Cauchy problem
	$$	\displaystyle \frac{du}{dt}+Au=0,\ \ t\ge0,\qquad
	u(0)=u_0.$$To this purpose,  consider in $L^1$ the nonlinear operator 
\begin{equation}\label{e3.1}
\begin{array}{rcl}
{A_0}y&=&-\Delta\beta(y)+{\rm div}(Db(y)y),\ \forall  y\in D({A_0}),\vsp
D(A_0)&=&\{y\in L^1,\ \beta(y)\in L^1_{\rm loc},-\Delta\beta(y)+{\rm div}(Db(y)y)\in L^1\},\end{array} \end{equation}where the differential operators $\Delta$ and div are taken in ${\mathcal{D}}'({\mathbb{R}}^d)$.
The main ingredient of the proof is the following lemma.
\begin{lemma}\label{l3.1} {We have}
	\begin{equation}
R(I+{{\lambda}} {A_0})=  L^1,\ \forall  {{\lambda}}>0,\label{e3.2}\end{equation}
and there is an operator $J_{{\lambda}} :L^1\to L^1$ such that
\begin{equation}
\begin{array}{ll}
J_{{\lambda}} u\in(I+{{\lambda}} A_0)^{-1} u,&\forall  u\in L^1,\vsp
|J_{{\lambda}} u-J_{{\lambda}} v|_1\le|u{-}v|_1,& \forall  u,v\in  L^1,{{\lambda}}>0.\end{array}\label{e3.3}
	\end{equation}
	\begin{equation}
	J_{{{\lambda}}_2} u=J_{{{\lambda}}_1}\left(\frac{{{\lambda}}_1}{{{\lambda}}_2}\,u
	+\left(1-\frac{{{\lambda}}_1}{{{\lambda}}_2}\right)J_{{{\lambda}}_2}u\right),\ \forall  \,0<{{\lambda}}_1,{{\lambda}}_2<{\infty}. \label{e3.3c}
		\end{equation}
	Moreover,
	\begin{equation}\label{e3.4a}
	\beta(J_{{\lambda}} u)\in L^q_{\rm loc},\  1<q<\frac d{d-1},
	\ \forall  u\in L^1,\end{equation}
	\begin{equation}\label{e3.3b}
\begin{array}{l}
|J_{{\lambda}} u|_{\infty}\le
(1+|({\rm div}\, D)^-+|D||^{\frac12}_{\infty}) |u|_{\infty},\ \forall  u\in L^1\cap L^{\infty},\vsp
\qquad0<{{\lambda}}<{{\lambda}}_0=
(|({\rm div}\,D)^-+|D||_\infty+(|({\rm div}\, D)^-+|D||^{\frac12}_{\infty})|b|_\infty)^{-1}, \end{array}\end{equation}
\begin{equation}\label{e3.6'}
|J_{{\lambda}}(u)|_{\infty}\le C_{{\lambda}}|u|_{\infty},\ \forall  \,u\in L^1\cap L^{\infty},\ \lambda>0,\mbox{ for some $C_{{\lambda}}\in(0,{\infty})$. }\end{equation}
\begin{equation}\label{e3.3bb}
J_{{\lambda}} ({\mathcal{P}}_0({\mathbb{R}}^d))\subset{\mathcal{P}}_0({\mathbb{R}}^d),  \ {{\lambda}}>0.\end{equation}
If \eqref{e2.4prim} holds, then
\begin{eqnarray}
&&{|J_{{\lambda}} g-g|_1\le C{{\lambda}}\|g\|_{W^{2,2}({\mathbb{R}}^d)},\ \forall  g\in C^{\infty}_0({\mathbb{R}}^d).}\label{e3.3a}
\end{eqnarray}
\end{lemma}

\noindent  Here $R(I+{{\lambda}} A_0)$ is the range of the operator $I+{{\lambda}} A_0$ and $(I+{{\lambda}} A_0)^{-1}:L^1\to D(A_0)$ (which, in general, might be multivalues) is a right  inverse of $A_0$. Before proving Lemma \ref{l3.1}, let us discuss some consequences.

Define the operator $A:D(A)\subset L^1\to L^1,$
\begin{eqnarray}
Au&=&A_0u,\ \ \forall  u\in D(A),\label{e3.5a}
\\
D(A)&=&\{J_{{\lambda}} v,\  v\in L^1\}.\label{e3.5aa}
\end{eqnarray}
	It is easily seen by \eqref{e3.3c} that
	\begin{equation}
	\label{e3.5aaa}
	D(A)\,=\,\{u=J_{{\lambda}} v,\ v\in L^1\},\ \forall \ 0<{{\lambda}}<{\infty}.
	\end{equation}
	By Lemma \ref{l3.1}, we have	
	\begin{lemma}\label{l3.2} \
		\begin{itemize}
			\item[{\rm(i)}] $J_{{\lambda}} $ coincides with the inverse $(I+{{\lambda}} A)^{-1}$ of $(I+{{\lambda}} A)$.
			\item[{\rm(ii)}] The operator $A$ is $m$-accretive in $L^1$, that is, $R(I+{{\lambda}} A)=L^1,$ $ \forall {{\lambda}}>0$, and
		\begin{equation}
		\label{e3.5aaaa}
	|(I+{{\lambda}} A)^{-1}u-(I+{{\lambda}} A)^{-1}v|_1\le|u-v|_1,\ \forall  u,v\in L^1,\ {{\lambda}}>0. \end{equation}
		\end{itemize}
		Moreover,  \eqref{e3.4a}, \eqref{e3.3b}, \eqref{e3.3bb} hold with $(I+{{\lambda}} A)^{-1}$ instead of $J_{{\lambda}} $. If \eqref{e2.4prim} holds, then $D(A)$ is dense in $L^1$.
	\end{lemma}
\noindent{\bf Proof of Lemma \ref{l3.2}.} (i) follows immediately by \eqref{e3.5a}, \eqref{e3.5aa},  \eqref{e3.5aaa}. Except for the density of $D(A)$ in $L^1$, all assertions of (ii) are immediate by the definition of $A$ and Lemma \ref{l3.1}. Now, assume that \eqref{e2.4prim} holds and   let   us prove that  $D(A)$ is dense in $L^1$ (that is, \mbox{$\overline{D(A)}=L^1$}).
By \eqref{e3.3a}, we have $C^{\infty}_0({\mathbb{R}}^d)\subset\overline{D(A)}$ and, since $C^{\infty}_0({\mathbb{R}}^d)$ is dense in $L^1$, the assertion follows.\bigskip

We recall that Lemma \ref{l3.2} implies via the Crandall and Liggett theorem (see \cite{1}, p.140) that, for each $u_0\!\in\!\overline{D(A)}{=}L^1$ and $T>0$, the Cauchy problem
\begin{equation}\label{e3.4}
\displaystyle \frac{du}{dt}+Au=0,\ t\in(0,T),\qquad
u(0)=u_0, \end{equation}has a unique mild solution $u\in C([0,T];L^1)$, that is,
 \begin{equation}\label{e3.5}
u(t)=\lim_{h\to0}u_h(t)\mbox{ in }L^1\mbox{ uniformly on }[0,T],\end{equation}where $u_h:[0,T]\to L^1$ is given by \eqref{e2.2}-\eqref{e2.4}, that is,
 \begin{equation}\label{e3.6}
 \begin{array}{l}
 u_h(t)=u^{i+1}_h,\ t\in(ih,(i+1)h],\vsp
 u^{i+1}_h+hAu^{i+1}_h=u^i_h,\ i=0,1,...,N-1;\ Nh=T,\vsp
 u^0_h=u_0.\end{array}\end{equation}
 In fact, the solution $u=u(t,u_0)$ given by \eqref{e3.5}, \eqref{e3.6} is given by the exponential formula 
\begin{equation}\label{e3.8b}
S(t)u_0=u(t,u_0)=\lim_{n\to{\infty}}
\left(I+\frac tn A\right)^{-n}u_0\mbox{ in }L^1,\ t\ge0,\end{equation}
where the convergence is uniform in $t$ on compact intervals $[0,T]$, and $S(t)$ is a semigroup of contractions on $L^1$, that is, 
$$\begin{array}{c}
S(t+s)=S(t)S(s)u_0,\ \forall  t,s\ge0,\quad
S(0)=I,\vsp
|S(t)u_0-S(t)\bar u_0|_1\le|u-\bar u_0|_1,\ \forall  u_0,\bar u_0\in L^1,\ t\ge0.\end{array}$$

\smallskip \noindent{\bf Proof of Lemma \ref{l3.1}.} We shall  follow the argument of \cite{6a} (Lemma 3.1). Namely, for $f\in L^1$ consider the equation\begin{equation}\label{e3.7}
u+{{\lambda}} {A_0}u= f \end{equation}or, equivalently,
 \begin{eqnarray}
& u-{{\lambda}} \Delta\beta(u)+{{\lambda}}\ {\rm div}(Db(u)u)=f\mbox{ in }{\mathcal{D}}'({\mathbb{R}}),\label{e3.8}\\[1mm]
& u\in L^1,\ \beta(u)\in L^1_{\rm loc},-\Delta\beta(u)+{\rm div}(Db(u)u)\in L^1.\label{e3.9}
\end{eqnarray} 
We shall assume first $f\in L^1\cap L^2$ and we approximate equation  \eqref{e3.8}~by
 \begin{equation}\label{e3.10}
u-{{\lambda}}\Delta\widetilde \beta_{\varepsilon}(u)+{{\lambda}}{\varepsilon}\widetilde \beta_{\varepsilon}(u)+{{\lambda}}\ {\rm div}(D_{\varepsilon}b_{\varepsilon}(u)u)=f,\end{equation}where $\widetilde \beta_{\varepsilon}(u)\equiv\beta_{\varepsilon}(u)+{\varepsilon} u$, and for ${\varepsilon}>0,\ r\in{\mathbb{R}}$,
\begin{equation}\label{e3.20'}\beta_{\varepsilon}(r)\equiv\frac1{\varepsilon}\ (r-(I+{\varepsilon}\beta)^{-1}r)=\beta((I+{\varepsilon}\beta)^{-1}r),\end{equation}
$$D_{\varepsilon}=\left\{\begin{array}{ll}
D,&\mbox{ if }|D|\in L^2\mbox{ and }{\rm div}\,D\in L^2+L^\infty,\\
\eta_{\varepsilon}D,&\mbox{ else}.\end{array}\right.$$
$$b_{\varepsilon}(r)=\left\{\begin{array}{ll}
b,&\mbox{\ if $b$ is a constant,}\vsp
\displaystyle \frac{(b*\rho_{\varepsilon})(r)}{1+{\varepsilon}|r|},&\mbox{ otherwise}.\end{array}\right.$$
Here $\rho_{\varepsilon}(r)=\frac1{\varepsilon}\ \rho\left(\frac r{\varepsilon}\right),\ \rho\in C^{\infty}_0({\mathbb{R}}),\ \rho\ge0$, is a standard modifier  and  by $I$ we denote the identity on ${\mathbb{R}}.$ Furthermore, $\eta_{\varepsilon}\in C^1_0(\mathbb{R}^d),$ $0\le\eta_{\varepsilon}\le1$, $|\nabla\eta_{\varepsilon}|\le1$ and $\eta_{\varepsilon}(x)=1$ if $|x|\le\frac1{\varepsilon}.$  Clearly, we have
\begin{eqnarray}\label{e3.23prim}
&|D_{\varepsilon}|\in L^\infty\cap L^2,\ |D_{\varepsilon}|\le|D|,\ \displaystyle\lim_{\varepsilon\to0} D_{\varepsilon}(x)=D(x),\ \mbox{ a.e. }x\in\mathbb{R}^d,\qquad\\ 
&\label{e3.23secund}
({\rm div}\,D_{\varepsilon})^-\le({\rm div}\,D)^-+{\bf1}_{|x|\ge\frac1\varepsilon}\,|D|. \end{eqnarray}

  We are going to show that, for ${\varepsilon}\to0$, the solution $\{u_{\varepsilon}\}$ to \eqref{e3.10} is convergent to a solution $u$ to \eqref{e3.8}. We can rewrite \eqref{e3.10} as
 \begin{equation}\label{e3.11}
({\varepsilon} I-\Delta)^{-1}u+{{\lambda}}\widetilde \beta_{\varepsilon}(u)
+{{\lambda}}({\varepsilon} I-\Delta)^{-1}{\rm div}(D_\varepsilon b_{\varepsilon}(u)u) =({\varepsilon} I-\Delta)^{-1}f.\end{equation}
We set
$$
\begin{array}{rcl}
F_{\varepsilon}(u)&=& ({\varepsilon} I-\Delta)^{-1}u,\  G^1_{\varepsilon}(u)={{\lambda}}\widetilde \beta_{\varepsilon}(u),\ u\in L^2,\vsp
G^2_{\varepsilon}(u)&=&{{\lambda}}({\varepsilon} I-\Delta)^{-1}({\rm div}(D_\varepsilon b_{\varepsilon}(u)u)),\ u\in L^2.\end{array}$$It is easily seen that $F_{\varepsilon}$ and $G^1_{\varepsilon}$ are accretive and continuous in $L^2.$

Since $r\mapsto b_{\varepsilon}(r)r$ is Lipschitz, we also have by assumptions (ii), (iii) that $G^2_{\varepsilon}$ is continuous in $L^2$ and
$$
\begin{array}{l}
\displaystyle \int_{{\mathbb{R}}^d}(G^2_{\varepsilon}(u)-G^2_{\varepsilon}(\bar u))(u-\bar u)dx\\
\qquad=-{{\lambda}}\displaystyle \int_{{\mathbb{R}}^d}D(b_{\varepsilon}(u)u-b_{\varepsilon}(\bar u)\bar u)
\cdot\nabla({\varepsilon} I-\Delta )^{-1}(u-\bar u)dx\vsp
\qquad\displaystyle \ge -C_{\varepsilon}{{\lambda}} |u-\bar u|_2|\nabla({\varepsilon} I-\Delta )^{-1}(u-\bar u)|_2,\ \forall  u,\bar u\in L^2({\mathbb{R}}^d).\end{array}$$ Moreover, we have
$$
\int_{{\mathbb{R}}^d}({\varepsilon} I-\Delta )^{-1}uu\,dx={\varepsilon}|({\varepsilon} I-\Delta )^{-1}u|^2_2+|\nabla({\varepsilon} I-\Delta )^{-1}u|^2_2,\ \forall  u\in L^2.$$
Hence, for $u^*=u-\bar u$, we have
$$\begin{array}{l}
(F_{\varepsilon}(u^*)+G^1_{\varepsilon}(u)-G^1_{\varepsilon}(\bar u)+G^2_{\varepsilon}(u)-G^2_{\varepsilon}(\bar u),u^*)_2\vsp
\qquad
\ge{{\lambda}}{\varepsilon}|u^*|^2_2+|\nabla({\varepsilon} I-\Delta )^{-1}u^*|^2_2+{\varepsilon}|({\varepsilon} I-\Delta )^{-1}u^*|^2_2\vsp
\qquad-C_{\varepsilon}{{\lambda}} |u^*|_2|\nabla({\varepsilon} I-\Delta )^{-1}u^*|_2.\end{array}$$This implies that $F_{\varepsilon}+G^1_{\varepsilon}+G^2_{\varepsilon}$ is  accretive and   coercive on $L^2$ for ${{\lambda}}<{{\lambda}}_{\varepsilon}$, where ${{\lambda}}_{\varepsilon}>0$ is sufficiently small. Since this operator is continuous and accretive, it follows that it is $m$-accretive and, therefore, surjective (because it is coercive). Hence, for each $f\in L^2\cap L^1$ and $0<{{\lambda}}<{{\lambda}}_{\varepsilon}$, equation \eqref{e3.11} has a unique solution $u_{\varepsilon}\in L^2$ with $\widetilde \beta_{\varepsilon}(u_{\varepsilon})\in H^1.$ Since $\widetilde \beta_{\varepsilon}$ has a Lipschitz inverse, we have that $u_{\varepsilon}\in H^1$, and hence  $b_{\varepsilon}(u_{\varepsilon})u_{\varepsilon}\in H^1\cap L^\infty$. 

We denote by $u_{\varepsilon}(f)\in H^1$ the solution to \eqref{e3.11} for $f\in L^2\cap L^1$ and we shall prove that
\begin{equation}\label{e3.12}
|u_{\varepsilon}(f_1)-u_{\varepsilon}(f_2)|_1\le|f_1-f_2|_1,\ \forall  f_1,f_2\in L^1\cap L^2.
\end{equation}To this purpose we set $u=u_{\varepsilon}(f_1)-u_{\varepsilon}(f_2),$ $f=f_1-f_2$.   By \eqref{e3.10}, we have, for $u_i=u_{\varepsilon}(f_i)$, $i=1,2,$
\begin{equation}\label{e3.13}
\begin{array}{r}
u-{{\lambda}}\Delta(\widetilde \beta_{\varepsilon}(u_1)-\widetilde \beta_{\varepsilon}(u_2))
+{\varepsilon}{{\lambda}}(\widetilde \beta_{\varepsilon}(u_1)-\widetilde \beta_{\varepsilon}(u_2))\vsp
\qquad+{{\lambda}}\,{\rm div}(D_{\varepsilon}(b_{\varepsilon}(u_1)u_1-b_{\varepsilon}(u_2)u_2))=f.\end{array}\end{equation}
We consider   the Lipschitzian function ${\mathcal{X}}_\delta:{\mathbb{R}}\to{\mathbb{R}},$
\begin{equation}\label{e3.18az}
{\mathcal{X}}_\delta(r)=\left\{\begin{array}{rl}
1&\mbox{ for }r\ge\delta,\vsp\displaystyle \frac r\delta&\mbox{ for }|r|<\delta,\vsp-1&\mbox{ for }r<-\delta,\end{array}\right.\end{equation}
where $\delta>0$. (We note that ${\mathcal{X}}_\delta(r)\to{\rm sign}\ r$ for $\delta\to0.$) We set
$$\Phi_{\varepsilon}={{\lambda}} \nabla(\widetilde \beta_{\varepsilon}(u_1)-\widetilde \beta_{\varepsilon}(u_2))
-{{\lambda}} D_{\varepsilon}(b_{\varepsilon}(u_1)u_1-b_{\varepsilon}(u_2)u_2)$$and rewrite \eqref{e3.13} as
$$u={\rm div}\ \Phi_{\varepsilon}-{\varepsilon}{{\lambda}}(\widetilde \beta_{\varepsilon}(u_1)-\widetilde \beta_{\varepsilon}(u_2))+f.$$
In particular, ${\rm div}\,\Phi_\varepsilon\in L^2$. We set $\Lambda_\delta={\mathcal{X}}_\delta(\widetilde \beta_{\varepsilon}(u_1)-\widetilde \beta_{\varepsilon}(u_2))$. Since $\Lambda_\delta\in H^1$, it follows that \mbox{$\Lambda_\delta\ {\rm div}\,\Phi_{\varepsilon}\in L^1$}, and so we have

\begin{equation}\label{e3.15}
\begin{array}{rl}
\displaystyle \int_{{\mathbb{R}}^d}u\Lambda_\delta dx=\!\!\!
&-\displaystyle \int_{{\mathbb{R}}^d}\Phi_{\varepsilon}\cdot\nabla
\Lambda_\delta dx\vsp&
-\,{\varepsilon}{{\lambda}}\displaystyle \int_{{\mathbb{R}}^d}(\widetilde \beta_{\varepsilon} (u_1)-\widetilde \beta_{\varepsilon} (u_2))\Lambda_\delta dx+\displaystyle \int_{{\mathbb{R}}^d}f\Lambda_\delta dx
 \vsp
=\!\!\!&-\displaystyle \int_{{\mathbb{R}}^d}\Phi_{\varepsilon}\cdot\nabla (\widetilde \beta_{\varepsilon}(u_1)-\widetilde \beta_{\varepsilon}(u_2)){\mathcal{X}}'_\delta(\widetilde \beta_{\varepsilon}(u_1)\vsp
&-\,\widetilde \beta_{\varepsilon}(u_2))dx-{\varepsilon}{{\lambda}}\displaystyle \int_{{\mathbb{R}}^d}(\widetilde \beta_{\varepsilon}(u_1)-\widetilde \beta_{\varepsilon}(u_2)){\mathcal{X}}_\delta(\widetilde \beta_{\varepsilon}(u_1)\vsp
&-\,\widetilde \beta_{\varepsilon}(u_2))dx
+\displaystyle \int_{{\mathbb{R}}^d}f\Lambda_\delta dx .\end{array}
\end{equation}
We have
$$\begin{array}{lcl}
I^1_{\delta}&=&\displaystyle \int_{{\mathbb{R}}^d}D_{\varepsilon}(b_{\varepsilon}(u_1)u_1-b_{\varepsilon}(u_2)u_2)\cdot\nabla\Lambda_\delta dx\vsp
&=&\displaystyle \frac1\delta
\int_{[|\widetilde \beta_{\varepsilon}(u_1)-\widetilde \beta_{\varepsilon}(u_2)|\le\delta]}D_{\varepsilon}(b_{\varepsilon}(u_1)u_1-b_{\varepsilon}(u_2)u_2)\cdot\nabla (\widetilde \beta_{\varepsilon}(u_1)-\widetilde \beta_{\varepsilon}(u_2))dx.\end{array}$$
Since  the inverse of $\widetilde \beta_{\varepsilon}$ is Lipschitz with constant $\frac1{\varepsilon}$,    we have
$$|b_{\varepsilon}(u_1)u_1-b_{\varepsilon}(u_2)u_2|\le C_{\varepsilon}|u_1-u_2|\le\frac{C_{\varepsilon}}{\varepsilon}\ |\widetilde \beta_{\varepsilon}(u_1)-\widetilde \beta_{\varepsilon}(u_2)|,$$and this yields
\begin{eqnarray}
&&\hspace*{-8mm}\displaystyle \lim_{\delta\to0} \frac1\delta\int_{[|\widetilde \beta_{\varepsilon}(u_1)-\widetilde \beta_{\varepsilon}(u_2)|\le\delta]}|D_{\varepsilon}(b_{\varepsilon}(u_1)u_1-b_{\varepsilon}(u_2)u_2)\cdot\nabla (\widetilde \beta_{\varepsilon}(u_1)-\widetilde \beta_{\varepsilon}(u_2))|dx\nonumber\\
&&\displaystyle \le \displaystyle \frac{C_{\varepsilon}}{\varepsilon}\ |D_{\varepsilon}|_2
\lim_{\delta\to0}\left(\int_{[|\widetilde \beta_{\varepsilon}(u_1)-\widetilde \beta_{\varepsilon}(u_2)|\le\delta]} |\nabla(\widetilde \beta_{\varepsilon}(u_1)-\widetilde \beta_{\varepsilon}(u_2))|^2dx\right)^{\frac12} \!\!=0,\qquad\quad\label{e3.24prim}\end{eqnarray}because $\nabla (\widetilde \beta_{\varepsilon}(u_1)-\widetilde \beta_{\varepsilon}(u_2))(x)=0$, a.e. on $[x\in{\mathbb{R}}^d;\widetilde \beta_{\varepsilon}(u_1(x))-\widetilde \beta_{\varepsilon}(u_2(x)){=}0]$.
This yields
 $ \lim\limits_{\delta\to0}I^1_{\delta}=0.$
Since ${\mathcal{X}}'_\delta\ge0$, a.e. on ${\mathbb{R}}$, we also have
$$\displaystyle \int_{{\mathbb{R}}^d}\nabla(\widetilde \beta_{\varepsilon}(u_1)-\widetilde \beta_{\varepsilon}(u_2))
\cdot\nabla (\widetilde \beta_{\varepsilon}(u_1)-\widetilde \beta_{\varepsilon}(u_2)){\mathcal{X}}'_\delta
(\widetilde \beta_{\varepsilon}(u_1)-\widetilde \beta_{\varepsilon}(u_2))\,dx\ge0.$$Then, by \eqref{e3.15},   we get
$$\lim_{\delta\to0}\int_{{\mathbb{R}}^d} u{\mathcal{X}}_\delta(\widetilde \beta_{\varepsilon}(u_1)-\widetilde \beta_{\varepsilon}(u_2))dx\le \int_{{\mathbb{R}}^d} |f|\,dx,\ \forall {\varepsilon}>0,$$and, since $u{\mathcal{X}}_\delta(\widetilde \beta_{\varepsilon}(u_1)-\widetilde \beta_{\varepsilon}(u_2))\ge0$ and ${\mathcal{X}}_\delta\to{\rm sign}$ as $\delta\to0$, by Fatou's lemma  this yields
\begin{equation}\label{e3.16}
|u|_1\le|f|_1,\end{equation}as claimed. We note that, since $f=0$ implies $u_{\varepsilon}(f)=0$, it follows by \eqref{e3.16} that $u_{\varepsilon}(f)\in L^1$,   $\forall  f\in L^1\cap L^2.$

Next, for $f$ arbitrary in $L^1$,   consider a sequence $\{f_n\}\subset L^2$ such that $f_n\to f$ strongly in $L^1.$
Let $\{u^n_{\varepsilon}\}\subset L^1\cap L^2$ be the corresponding solutions to \eqref{e3.10} for $0<{{\lambda}}<{{\lambda}}_{\varepsilon}$.
Taking into account  \eqref{e3.16}, we obtain by the above equation that
 $|u^n_{\varepsilon}-u^m_{\varepsilon}|_1\le|f_n-f_m|_1,\ \forall  n,m\in{\mathcal{N}}.$ Hence, for $n\to{\infty}$, we have
$u^n_{\varepsilon}\to u_{\varepsilon}(f)\mbox{ in }L^1$. Define the operator
\begin{equation}\label{e3.25b}
\begin{array}{rcl}
A_{\varepsilon} u&=&-\Delta\widetilde \beta_{\varepsilon}(u)+{\varepsilon}\widetilde \beta_{\varepsilon}(u)+{\rm div}(D_{\varepsilon}b_{\varepsilon}(u)u)\vsp
D(A_{\varepsilon})&=&\{u\in L^1;-\Delta\widetilde \beta_{\varepsilon}(u)+{\varepsilon}\widetilde \beta_\varepsilon(u)+{\rm div}(D_{\varepsilon}b_{\varepsilon}(u)u)\in L^1\}.\end{array}\end{equation} It is obvious that $(A_{\varepsilon},D(A_{\varepsilon}))$ is closed on $L^1.$ Therefore,  $u_{\varepsilon}(f)\in D(A_{\varepsilon})$ and
\begin{equation}\label{e3.17}
u_{\varepsilon}(f)+{{\lambda}} A_{\varepsilon} u_{\varepsilon} (f)=f,
\end{equation}for ${{\lambda}}<{{\lambda}}_{\varepsilon}$. We also have
\begin{equation}\label{e3.18}
|u_{\varepsilon}(f_1)-u_{\varepsilon}(f_2)|_1\le|f_1-f_2|_1,\ \forall  f_1,f_2\in L^1.\end{equation}
Then, by Proposition 3.3 in \cite{1}, p.~99, it follows that  \mbox{$R(1+{{\lambda}} A_{\varepsilon})=L^1,$} \mbox{$ \forall {{\lambda}}>0,$} and also \eqref{e3.17}, \eqref{e3.18} extend to all ${{\lambda}}>0$ (see Proposition 3.1 in \cite{1}).

Moreover, if $u_{\varepsilon}=u_{\varepsilon}({{\lambda}},f)$ is our solution to \eqref{e3.10}, we have, for all \mbox{$0<{{\lambda}}_1,{{\lambda}}_2<{\infty}$} and $f\in L^1\cap L^2$, by definition
\begin{equation}
\label{e3.29a}
u_{\varepsilon}({{\lambda}}_2,f)=u_{\varepsilon}\left({{\lambda}}_1,\frac{{{\lambda}}_1}{{{\lambda}}_2}\,f+\left(1-\frac{{{\lambda}}_1}{{{\lambda}}_2}\right)u_{\varepsilon}({{\lambda}}_2,f)\right).	\end{equation}
If $f\in L^1\cap L^{\infty}$,  we have
\begin{equation}\label{e3.20a}
|u_{\varepsilon}(f)|_{\infty}\le
(1+|({\rm div}\, D)^-+|D||^{\frac12}_{\infty})
|f|_{\infty} ,\  \ 0<{{\lambda}}<{{\lambda}}_0.  \end{equation}
Indeed, by \eqref{e3.10}, we see that, for $M_{\varepsilon}=|({\rm div }\,D_{\varepsilon})^-|^{\frac12}_\infty|f|_\infty$,   $u_{\varepsilon}=u_{\varepsilon}(f)$, $b^*_{\varepsilon}(r)=b_{\varepsilon}(r)r$ and ${{\lambda}}<{{\lambda}}_0$, where $\lambda_0$ is as in \eqref{e3.3b},
$$\begin{array}{l}
(u_{\varepsilon}-|f|_{\infty}-M_{\varepsilon})
-{{\lambda}}\Delta(\widetilde \beta_{\varepsilon}(u_{\varepsilon})-\widetilde \beta_{\varepsilon}(|f|_{\infty}+M_{\varepsilon}))\vsp
\qquad +{{\lambda}}{\varepsilon}(\widetilde \beta_{\varepsilon}(u_{\varepsilon})-\widetilde \beta_{\varepsilon}(|f|_{\infty}+M_{\varepsilon}))
+{{\lambda}}\,{\rm div}(D_{\varepsilon}(b^*_{\varepsilon}(u_{\varepsilon})-b^*_{\varepsilon}(|f|_{\infty}+M_{\varepsilon})))\vsp
\qquad\le f-|f|_{\infty}-M_{\varepsilon}
-{{\lambda}} b^*_{\varepsilon}(M_{\varepsilon}+|f|_{\infty}){\rm div}\,D_{\varepsilon}\le0.
\end{array}$$
Multiplying the above equation by ${\mathcal{X}}_\delta((u_{\varepsilon}-
(|f|_{\infty}+M\red{_{\varepsilon}}))^+)$ and integrating over ${\mathbb{R}}^d$, we get as above, for $\delta\to0$,
 $|(u_{\varepsilon}-|f|_{\infty}-M\red{_{\varepsilon}})^+|_1\le0$
and, therefore,  by \eqref{e3.23secund} 
$$u_{\varepsilon}\le
 (1+|({\rm div}\, D)^-+|D||^{\frac12}_{\infty})|f|_{\infty},\   \mbox{ a.e. in }{\mathbb{R}}^d.$$Similarly, one gets that
$$u_{\varepsilon}\ge-(1+|({\rm div}\, D)^-+|D||^{\frac12}_{\infty})|f|_{\infty},\  \mbox{ a.e. in }{\mathbb{R}}^d,$$and so \eqref{e3.20a} follows, which in turn by \eqref{e3.29a} implies that, for some $C_\lambda>0$,
	\begin{equation}\label{e3.36prim}
	|u_\varepsilon(\lambda,f)|_\infty\le C_\lambda|f|_\infty\mbox{\ \ for all }\varepsilon>0,\ f\in L^1\cap L^\infty.\end{equation}

Now, we are going to let ${\varepsilon}\to0$ in \eqref{e3.17}.
To this end, we need some estimates on $u_{\varepsilon}$.
By \eqref{e3.10} (or \eqref{e3.17}), we have
\begin{equation}\label{e3.18a}
\Delta\widetilde \beta_{\varepsilon}(u_{\varepsilon})={{\lambda}}^{-1}(u_{\varepsilon}-f)+{\rm div}(D_{\varepsilon}b_{\varepsilon}(u_{\varepsilon})u_{\varepsilon})+{\varepsilon}\widetilde \beta_{\varepsilon}(u_{\varepsilon}).\end{equation}We shall first consider the case $d\ge3.$ This yields, if $f\in L^1\cap L^\infty$,\begin{equation}
\label{e3.19}
\widetilde \beta_{\varepsilon}(u_{\varepsilon})=\frac1{{\lambda}}\ E*(-u_{\varepsilon}
-{\varepsilon}{{\lambda}}\widetilde \beta_{\varepsilon}(u_{\varepsilon})+f)
+{\rm div} (E*(D_{\varepsilon}b_{\varepsilon}(u_{\varepsilon})u_{\varepsilon})),\mbox{ a.e. in }{\mathbb{R}}^d,
\end{equation}where $E$ is the fundamental solution to the Laplace operator (see Section \ref{s1}).
Then, recalling \eqref{e1.4}, \eqref{e1.5} and \eqref{e3.16}, we get, for ${\varepsilon}\in(0,1)$,
\begin{equation}
\label{e3.20}
\begin{array}{ll}
\|\widetilde \beta_{\varepsilon}(u_{\varepsilon}) -{\rm div}
(E*(D_{\varepsilon}b_{\varepsilon}(u_{\varepsilon})u_{\varepsilon}))\|_{M^{\frac d{d-2}}}
\\\qquad\quad\le\displaystyle \frac1{{\lambda}}\,\|E\|_{M^{\frac d{d-2}}}|u_{\varepsilon}+{{\lambda}}{\varepsilon}\widetilde \beta_{\varepsilon}(u_{\varepsilon})- f|_1
\le\displaystyle \frac{2+3{{\lambda}}}{{\lambda}}\,\|E\|_{M^{\frac d{d-2}}}|f|_1, \end{array}
\end{equation}$\forall {{\lambda}}>0,$ because ${\varepsilon}|\widetilde \beta_{\varepsilon}(u)|\le(2+{\varepsilon}^2)|u|.$
Taking into account that, for $d\ge3$,
$$
M^{\frac d{d-2}}\subset L^p_{\rm loc},\ \forall  p\in\mbox{$\left(1,\frac d{d-2}\right)$},\ \
M^{\frac d{d-1}} \subset L^p_{\rm loc},\ \forall  p\in\mbox{$\left(1,\frac d{d-1}\right)$} $$
we see by \eqref{e3.20} that, for $1<p<\frac d{d-2},$ and, for all compacts $K\subset{\mathbb{R}}^d$, we~have
$$\|\widetilde \beta_{\varepsilon}(u_{\varepsilon})-{\rm div} (E*(D_{\varepsilon}b_{\varepsilon}(u_{\varepsilon})u_{\varepsilon}))\|_{L^p(K)}
\le\frac{1+{{\lambda}}}{{\lambda}}\ C_K|f|_1,\ \forall {{\lambda}}>0,$$and so, by \eqref{e1.6}, we have, for $1<q<\frac d{d-1},$ $f\in L^1\cap L^\infty$,
\begin{equation}\label{e3.21}
\begin{array}{ll}
\|\widetilde \beta_{\varepsilon}(u_{\varepsilon})\|_{L^q(K)}\!\!\!
&\le C_K \left(\|{\rm div} (E*(D_{\varepsilon}b_{\varepsilon}(u_{\varepsilon})u_{\varepsilon}))\|_{L^q(K)}+\displaystyle \frac{1+{{\lambda}}}{{\lambda}}|f|_1\right)\vsp
&\le C_K \left(\|{\rm div} (E*(D_{\varepsilon}b_{\varepsilon}(u_{\varepsilon})u_{\varepsilon}))\|_{M^{\frac d{d-1}}}+\displaystyle \frac{1+{{\lambda}}}{{\lambda}}|f|_1\right)\vsp
&\le C_K \left(|u_{\varepsilon}|_1+\displaystyle \frac{1+{{\lambda}}}{{\lambda}}|f|_1\right)  \le\displaystyle  C_K\left(\displaystyle \frac{1+{{\lambda}}}{{\lambda}}\right) |f|_1,  \end{array}\end{equation}$\forall {{\lambda}}>0,$
for any compact subset $K\subset{\mathbb{R}}^d$, where the constant $C_K$ is independent of $f$ and changes from line to line and we used that $|b_{\varepsilon}|_{\infty}\le|b|_{\infty}$,  $|D_{\varepsilon}|\le|D|_\infty$. 

We assume first that $f\in L^1\cap L^{\infty}$and $0<{{\lambda}}<{{\lambda}}_0$. Then, by \eqref{e3.20a}, we have, along a subsequence $\{{\varepsilon}\}\to0$,
 \begin{equation}\label{e3.24}
u_{\varepsilon}\to u\mbox{ weak-star in $L^{\infty}$, whence weakly in }L^1_{\rm loc}.\end{equation}
Therefore,  $\{\widetilde \beta_{\varepsilon}(u_{\varepsilon})\}$ is bounded in $L^q_{\rm loc}$ (as a matter of fact, since $|\widetilde \beta_{\varepsilon}(r)|\le C_M(1+|r|),$ $|r|\le M$, ${\varepsilon}\in(0,1)$, see \eqref{e3.41prim} below, it is bounded in $L^{\infty}$). Hence,    by selecting  a further subsequence $\{{\varepsilon}\}\to0$, we  also have
 \begin{equation}\label{e3.25a}
\widetilde \beta_{\varepsilon}(u_{\varepsilon})\to \eta \mbox{ weakly in $L^q_{\rm loc},$ $q\in \left(1,\frac d{d-1}\right).$}\end{equation}
Now, we consider two cases:

\medskip \noindent{\bf Case 1.} $b$ is constant (hence $b_{\varepsilon}\equiv b$).\medskip

\noindent By
 \eqref{e3.10},  \eqref{e3.23prim}, \eqref{e3.24}, \eqref{e3.25a}, we have
 \begin{equation}\label{e3.26a}
 u-{{\lambda}}\Delta\eta+{{\lambda}}\ {\rm div}(Dbu)=f\mbox{ in }{\mathcal{D}}'({\mathbb{R}}^d).\end{equation}
It remains to be shown that
\begin{equation}\label{e3.25}\eta(x)=\beta(u(x)),\mbox{ a.e. }x\in{\mathbb{R}}^d.\end{equation}
For this purpose, we shall prove first, via the Riesz-Kolmogorov compactness theorem, that $\{\widetilde \beta_{\varepsilon}(u_{\varepsilon})\}_{{\varepsilon}>0}$ is compact in $L^q_{\rm loc}.$ We set
$$v_{\varepsilon}=\widetilde \beta_{\varepsilon}(u_{\varepsilon}),\ v^\nu_{\varepsilon}(x)=v_{\varepsilon}(x+\nu)-v_{\varepsilon}(x),\ \forall  x,\nu\in{\mathbb{R}}^d.$$By \eqref{e3.19}, we have
$$v^\nu_{\varepsilon}=\frac1{{\lambda}}\ E^\nu*(-u_{\varepsilon}-{\varepsilon}{{\lambda}}\widetilde \beta_{\varepsilon}(u_{\varepsilon})+f)-{\rm div} (E^\nu*(D_{\varepsilon}bu_{\varepsilon})),\ E^\nu(x)\equiv E(x+\nu)-E(x),$$and, by \eqref{e1.5}, \eqref{e1.6},  this yields, for any compact $K\subset{\mathbb{R}}^d$,  $1<q<\frac d{d-1},$ $q<p<\frac d{d-2},$ and all ${\varepsilon}\in(0,1)$,
$$\begin{array}{ll}
\|v^\nu_{\varepsilon}\|_{L^q(K)}\!\!\!
&\le\displaystyle \frac{C_K}{{\lambda}}\,\|E^\nu*(-u_{\varepsilon}-{\varepsilon}{{\lambda}}\widetilde \beta_{\varepsilon}(u_{\varepsilon})+f)\|_{L^p(K)}
\vsp&\qquad+\|{\rm div} (E^\nu*(D_{\varepsilon}bu_{\varepsilon}))\|_{L^q(K)} \vsp
&\le\displaystyle \frac{C_K}{{\lambda}}\,\|E^\nu*(-u_{\varepsilon}-{\varepsilon}{{\lambda}}\widetilde \beta_{\varepsilon}(u_{\varepsilon})+f)\|_{M^{\frac d{d-2}}}\vsp&\qquad
+C_K\|{\rm div} (E^\nu*(D_{\varepsilon}bu_{\varepsilon}))\|_{M^{\frac d{d-1}}}\vsp
&\le\displaystyle \frac{C_K}{{\lambda}}\,\|E^\nu\|_{M^{\frac d{d-2}}}(|u_{\varepsilon}|_1+|f|_1)
+C_K\|\nabla E^\nu\|_{M^{\frac d{d-1}}}|u_{\varepsilon}|_1\vsp
&\le C_K\left(1+\displaystyle \frac1{{\lambda}}\right)(\|E^\nu\|_{M^{\frac d{d-2}}}+\|\nabla E^\nu\|_{M^{\frac d{d-1}}})|f|_1,\end{array}$$ where we used that $|D_\varepsilon|_\infty\le|D|_\infty$. 
On the other hand, we have
$$\lim_{\nu\to0}(\|E^\nu\|_{M^{\frac d{d-2}}}+
\|\nabla E^\nu\|_{M^{\frac d{d-1}}})=0.$$
(This continuity property follows as in the case of $L^p$-spaces.)
Therefore, $\{\widetilde \beta_{\varepsilon}(u_{\varepsilon});\ {\varepsilon}\in(0,1)\}$ is compact in each space $L^q(K)$, $K$ compact subset of ${\mathbb{R}}^d$, where \mbox{$1<q<\frac d{d-1}.$} 
We also note that ${\varepsilon}|u_{\varepsilon}|_1\to0$ as ${\varepsilon}\to0$. Hence, on a subsequence $\{{\varepsilon}\}\to0$,
\begin{equation}\label{e3.26}
\beta_{\varepsilon}(u_{\varepsilon})\to\eta\mbox{ strongly in $L^q_{\rm loc}$}.
\end{equation}
Since $\{u_{\varepsilon};\ {\varepsilon}\in(0,1)\}$ are bounded in $L^{\infty}$, we have
$$\lim_{{\varepsilon}\to0}|\beta_{\varepsilon}(u_{\varepsilon})-\beta(u_{\varepsilon})|_{\infty}=0,$$so, by \eqref{e3.26},
\begin{equation}
\label{e3.37a}
\beta(u_{\varepsilon})\to\eta\mbox{ strongly in }L^q_{\rm loc}.\end{equation}	
Recalling that $u_{\varepsilon}\to u$ weak-star in $L^{\infty}$ and that the map $u\to\beta(u)$ is maximal monotone in each dual pair $(L^q(K),L^{q'}(K))$, hence weakly-strongly closed, we get by \eqref{e3.37a}  that \eqref{e3.25} holds. Hence, by \eqref{e3.18}, \eqref{e3.24},  for   $u=u({{\lambda}},f)$,  we~have
\begin{equation}\label{e3.27}
|u({{\lambda}},f)-u({{\lambda}},g)|_1\le|f-g|_1,\ \forall {{\lambda}}>0,\ f,g\in L^1\cap L^{\infty}.
\end{equation}
(Indeed, first only for $0<{{\lambda}}<{{\lambda}}_0$,   but then by Proposition 3.1 in \cite{1} for all ${{\lambda}}>0.$)
Hence, $u(\lambda,f)\in L^1$ and so, by \eqref{e3.26a}, $u(\lambda,f)\in D(A_0)$ and satisfies \eqref{e3.7} for all $f\in L^1\cap L^\infty$.

\medskip \noindent{\bf Case 2.} Let $\beta$ be strictly increasing.\smallskip

\noindent  Multiplying \eqref{e3.17} by $u_{\varepsilon}=u_{\varepsilon}(f)$ and integrating over ${\mathbb{R}}^d$, since $u_{\varepsilon},b_{\varepsilon}(u_{\varepsilon})u_{\varepsilon}\in H^1\cap L^1\cap L^{\infty},$ $\widetilde \beta_{\varepsilon}(u_{\varepsilon})\in H^2$, we obtain
\begin{equation}
\label{e3.40aa}
\begin{array}{l}
|u_{\varepsilon}|^2_2+{{\lambda}}\displaystyle \int_{{\mathbb{R}}^d}
\beta'_{\varepsilon}(u_{\varepsilon})|\nabla u_{\varepsilon}|^2dx\vsp\qquad\
\le{{\lambda}}\displaystyle \int_{{\mathbb{R}}^d}(\nabla u_{\varepsilon}\cdot D_{\varepsilon})b_{\varepsilon}(u_{\varepsilon})u_{\varepsilon} dx
+\displaystyle \frac12\ |u_{\varepsilon}|^2_2+\frac12\ |f|^2_2.
\end{array}\end{equation}
Defining
$$\psi(r)=\int^r_0b_{\varepsilon}(s)s\,ds,\ r\in{\mathbb{R}},$$we see that $\psi\ge0$, hence the first integral on the right hand side of \eqref{e3.40aa} is equal to
\begin{equation}
\label{e3.40b}
-\int_{\mathbb{R}}{\rm div}\,D_{\varepsilon}\psi(u_{\varepsilon})dx\le C<\infty,\end{equation}where (see \eqref{e3.23secund})
	$$C:=\frac12\,|b|_\infty|({\rm div}\,D)^-+|D||_\infty\sup_{\varepsilon\in(0,1)}(|u_\varepsilon|_\infty|u_\varepsilon|_1).$$
	$C$ is by \eqref{e3.18} and \eqref{e3.36prim} indeed finite, since $f\in L^1\cap L^\infty.$

Define $g_{\varepsilon}(r)=(I+{\varepsilon}\beta)^{-1}(r),$ $r\in{\mathcal{N}}$, and
$$a(r)=\int^r_0\frac{\beta'(s)}{1+\beta'(s)}\ ds,\ \ r\in{\mathbb{R}}.$$Since
$$\beta'_{\varepsilon}(r)\ge\frac{\beta'(g_{\varepsilon}(r))}{1+\beta'(g_{\varepsilon}(r))}
\ge\frac{\beta'(g_{\varepsilon}(r))}{1+\beta'(g_{\varepsilon}(r))}\ (g'_{\varepsilon}(r))^2,\ r\in{\mathbb{R}},$$and thus
$$\beta'_{\varepsilon}(u_{\varepsilon})|\nabla u_{\varepsilon}|^2\ge|\nabla a(g_{\varepsilon}(u_{\varepsilon}))|^2,$$we obtain from \eqref{e3.40aa}, \eqref{e3.40b}
\begin{equation}
\label{e3.40c}
|u_{\varepsilon}|^2_2+2{{\lambda}}\int_{{\mathbb{R}}^d}|\nabla a(g_{\varepsilon}(u_{\varepsilon}))|^2dx\le|f|^2_2.\end{equation}
Since $|a(r)|\le|r|$ and $|g_{\varepsilon}(r)|\le |r|$,
 $r\in{\mathbb{R}}$, this implies that $a(g_{\varepsilon}(u_{\varepsilon}))$, ${\varepsilon}>0$, is bounded in $H^1$, hence compact in $L^2_{\rm loc}$, so along a subsequence ${\varepsilon}\to0$
$$a(g_{\varepsilon}(u_{\varepsilon}))\to v\mbox{\ \ in $L^2_{\rm loc}$ and a.e.}$$ Since $a$ is strictly increasing and continuous, and thus so is its inverse function $a^{-1}$, it follows that
 $g_{\varepsilon}(u_{\varepsilon})\to a^{-1}(v),\mbox{ a.e. },$ and so, as ${\varepsilon}\to0,$
	$$u_{\varepsilon}=g_{\varepsilon}(u_{\varepsilon})+{\varepsilon}\beta(g_{\varepsilon}(u_{\varepsilon}))\to a^{-1}(v),\mbox{ a.e. on }{\mathbb{R}}^d.$$Therefore, by \eqref{e3.20a} we have $u\in L^{\infty}$ and by \eqref{e3.24}
\begin{equation}
\label{e3.40d}
u_{\varepsilon}\ \st{{\varepsilon}\to0}\longrightarrow\ u\mbox{\ \ in }L^p_{\rm loc},\ \ \forall  p\in[1,{\infty}).\end{equation}
By Fatou's lemma and \eqref{e3.18}, it follows that $u\in L^1$. Furthermore, obviously both $\widetilde \beta_{\varepsilon}$ and $b_{\varepsilon}$, ${\varepsilon}\in(0,1)$, are locally equicontinuous. Hence,
\begin{equation}\label{e3.50'}
\widetilde \beta_{\varepsilon}(u_{\varepsilon})\to\beta(u),\ \ b_{\varepsilon}(u_{\varepsilon})u_{\varepsilon}\to b(u)u,\mbox{ a.e. on }{\mathbb{R}}^d,\end{equation}as ${\varepsilon}\to0,$ since $\sup\limits_{{\varepsilon}>0}|u_{\varepsilon}|_{\infty}<{\infty}$ by \eqref{e3.20a}. We note that, since $\beta$ is locally Lipschitz and $\beta(0)=0$, we have, for $M>0$, $C_M=\sup\limits_{|r|\le M}\beta'(r)$,
\begin{equation}
\label{e3.50''}
|\beta_{\varepsilon}(r)|\le C_M|(I+{\varepsilon}\beta)^{-1}r|\le C_M|r|,\ r\in[-M,M].
\end{equation}Hence, by \eqref{e3.20a}, \eqref{e3.36prim}  and because $|b_{\varepsilon}|\le|b|_{\infty}$, \eqref{e3.50''} implies that both convergences in \eqref{e3.50'} also hold  in $L^p_{\rm loc}$, $p\in[1,{\infty})$. Hence we can pass to the limit in \eqref{e3.17}, \eqref{e3.18} to conclude that $u$ satisfies \eqref{e3.8}, which in turn implies that $u\in D(A),$ because $\beta(u)\in L^p_{\rm loc},$ $p\in[1,{\infty}),$ and that \eqref{e3.27} also holds in Case 2.

All what follows now holds in both Cases 1 and 2. 

We define $J_{{\lambda}} :L^1\cap L^{\infty}\to L^1$,
$$J_{{\lambda}} (f)=u({{\lambda}},f),\ \ \forall \,f\in L^1\cap L^{\infty},\ {{\lambda}}>0.$$Then,   \eqref{e3.3c} follows by definition.   We note that, by 
 \eqref{e3.21}, it also follows
\begin{equation}\label{e3.35a}
|\beta(u({{\lambda}},f))|_{L^q(K)}
\le C_K\left(1+\frac1{{\lambda}}\right)|f|_1,\ \forall  f\in L^1\cap L^{\infty},\ {{\lambda}}>0.\end{equation}Now, let $f\in L^1$ and $\{f_n\}\subset L^1\cap L^{\infty}$ be such that $f_n\to f$  in $L^1$. If $u_n=u({{\lambda}},f_n)$, we see by \eqref{e3.27} that $u_n\to u$ in $L^1$ and, by \eqref{e3.35a}, $\{\beta(u_n)\}$ is bounded in $L^q(K)$ for each $K\subset {\mathbb{R}}^d$, $q\in\left(0,\frac d{d-1}\right)$. Hence, by the generalized Lebesgue  convergence theorem and the continuity of $\beta$,  $\beta(u_n)\to\beta(u)$   in $L^q(K)$, and so $u\in D(A_0)$, $A_0u_n\to A_0u$ in $L^1$, and
\begin{equation}
\label{e3.41prim}
u+{{\lambda}} A_0u=f,
\end{equation}Then, we may extend \eqref{e3.27} and \eqref{e3.35a} to all of $f\in L^1$ and hence \eqref{e3.3} and \eqref{e3.4a} hold.

By \eqref{e3.20a}, \eqref{e3.24},  it follows that \eqref{e3.3b} holds, which by \eqref{e3.3c} implies \eqref{e3.6'}. 
Moreover,if $f\ge0$ on ${\mathbb{R}}^d$,  by \eqref{e3.10} it is easily seen that $u_{\varepsilon}\ge0$ on ${\mathbb{R}}^d$, and so, by \eqref{e3.24} we infer that $u\ge0$. Moreover, we have
\begin{equation}\label{e3.31b}
\int_{{\mathbb{R}}^d}u\,dx=\int_{{\mathbb{R}}^d}f\,dx.
\end{equation}
By \eqref{e3.41prim},   we have
\begin{equation}\label{e3.31aa}
 \displaystyle \int_{{\mathbb{R}}^d} u{\varphi} dx-{{\lambda}}\int_{{\mathbb{R}}^d} \beta(u)\Delta{\varphi} dx+{{\lambda}}\int_{{\mathbb{R}}^d} ub(u)   D\cdot\nabla{\varphi} dx \displaystyle
 =\displaystyle \int_{{\mathbb{R}}^d}\!\!f{\varphi}\,dx,
\end{equation}
$\forall {\varphi}\in C^{\infty}_0.$ Now, we choose in \eqref{e3.31aa} ${\varphi}={\varphi}_\nu\in C^{\infty}_0$, where ${\varphi}_\nu\to1$ on   ${\mathbb{R}}^d$, $0\le{\varphi}_\nu\le1$, and $$\mbox{$|\Delta{\varphi}_\nu|_{\infty}+|\nabla{\varphi}_\nu|_{\infty}\to0$ as $\nu\to0.$}$$
(Such an example is ${\varphi}_\nu(x)=\exp\left(-\frac{\nu|x|^2}{1-\nu|x|^2}\right).$)
This implies \eqref{e3.31b}, and so  \eqref{e3.3bb} holds.

It remains to prove \eqref{e3.3a}. Let $g\in C^{\infty}_0({\mathbb{R}}^d)$ be arbitrary but fixed. Coming back to equation \eqref{e3.17}, we can write
	$$g-{{\lambda}}\Delta\widetilde \beta_{\varepsilon}(g)+
	{{\lambda}}{\varepsilon}\widetilde \beta_{\varepsilon}(g)+{{\lambda}}\ {\rm div}(D_{\varepsilon}b_{\varepsilon}(g)g)=g+{{\lambda}} A_{\varepsilon}(g),$$and so \eqref{e3.18} yields that, for some $C\in(0,\infty)$,
	$$|u_{\varepsilon}(g)-g|_1\le{{\lambda}}|A_{\varepsilon}(g)|_1\le C{{\lambda}}(\|g\|_{W^{2,2}({\mathbb{R}}^d)}+|g|_\infty),\ \forall\varepsilon\in(0,1),$$because $\beta\in C^2$, $b\in C^1$, and hence $\widetilde \beta_{\varepsilon}$, $(\widetilde \beta_{\varepsilon})'$, $(\widetilde \beta_{\varepsilon})''$ and $b'_\varepsilon$ are locally uniformly bounded in ${\varepsilon}\in(0,1)$ and since $|D|_\infty<\infty$, ${\rm div}\,D\in L^1_{\rm loc}$ and $|b_{\varepsilon}|\le|b|_{\infty}$. This, together with \eqref{e3.24}, implies \eqref{e3.3a}, as claimed.
  This completes the proof of Lemma \ref{l3.1} for $d\ge3.$

  Now, we shall sketch the proof of the case $d=1$ only in the case that $b\equiv1$.
By \eqref{e3.18a}, we have
$$(\widetilde \beta_{\varepsilon}(u_{\varepsilon}))'(x)=D(x)u_{\varepsilon}(x)+{{\lambda}}^{-1}\int^x_{-{\infty}}(u_{\varepsilon}(y)-f(y))dy,\ \forall  x\in{\mathbb{R}},$$
\begin{equation}
\label{e3.31a}
\begin{array}{ll}
\widetilde \beta_{\varepsilon}(u_{\varepsilon}(x))\!\!\!&=\widetilde \beta_{\varepsilon}(u_{\varepsilon}(x_0))+\displaystyle \int^x_{x_0}D(y)u_{\varepsilon}(y)dy\\
&+{{\lambda}}^{-1}\displaystyle \int^x_{x_0}ds
\int^s_{-{\infty}}(u_{\varepsilon}(y)-f(y))dy,\ \forall  x,x_0\in{\mathbb{R}}.\end{array}
\end{equation}
Taking into account that $\{u_{\varepsilon}\}$ is bounded in $L^1$, we may choose $x_0$ independent of ${\varepsilon}$ such that $\{u_{\varepsilon}(x_0)\}$ is bounded. This implies that $\{\widetilde \beta_{\varepsilon}(u_{\varepsilon})\}_{{\varepsilon}>0}$ is bounded in $L^{\infty}_{\rm loc}$ and so estimate \eqref{e3.21}  follows. Hence, it follows as above that \eqref{e3.24}-\eqref{e3.26a} follow too. We shall prove that $\{\widetilde \beta_{\varepsilon}(u_{\varepsilon})\}_{{\varepsilon}>0}$ is compact in $L^1_{\rm loc}$. If $v_{\varepsilon}= \beta_{\varepsilon}(u_{\varepsilon})$ and $v^\nu_{\varepsilon}=v_{\varepsilon}(x+\nu)-v_{\varepsilon}(x)$, we get by \eqref{e3.31a}~that
$$v^\nu_{\varepsilon}(x)=\int^{x+\nu}_xD(y)u_{\varepsilon}(y)dy
+{{\lambda}}^{-1}\displaystyle \int^{x+\nu}_xds\int^s_{-{\infty}}(u_{\varepsilon}(y)-f(y))dy$$and, therefore, $\lim\limits_{\nu\to0}|v^\nu_{\varepsilon}|_{L^1_{\rm loc}}=0$ uniformly with respect to ${\varepsilon}$.  Hence, $\{v_{\varepsilon}\}$ is compact in $L^1_{\rm loc}$ and so \eqref{e3.25} follows by \eqref{e3.26}. As regard \eqref{e3.27}-\eqref{e3.31b}, and so all  conclusion of Lemma \ref{l3.1}, it follows as in the previous case.

Later, we shall also need the following convergence result for the solution $u_{\varepsilon}$ to the approximating equation \eqref{e3.10}.

\begin{lemma}\label{l3.2a} Assume that $\beta$ is strictly increasing. Then, we have, for ${\varepsilon}\to0$,
\begin{equation}\label{e3.52a}
u_{\varepsilon}\to u=J_{{\lambda}} f\mbox{ in }L^1,\ \forall \ f\in L^1.\end{equation}
\end{lemma}

\noindent{\bf Proof.} We shall proceed as in the proof of Lemma 3.2 in \cite{6a}. It suffices to prove this for $f$ in a dense subset of $L^1$. In Case 2 of the proof of Lemma \ref{l3.1}, by \eqref{e3.40d} we have  $u_{\varepsilon}\to u=J_{{\lambda}} f$ strongly in $L^1_{\rm loc}$. But this also follows in Case 1, because by \eqref{e3.25} and \eqref{e3.37a} we have that $\beta(u_{\varepsilon})\to\beta(u)$ in $L^1_{\rm loc}$. So, by our additional assumption that  $\beta$ is strictly increasing, which entails that its inverse $\beta^{-1}$ is continuous, we have that also $u_{\varepsilon}\to u$, a.e. (along a subsequence), hence $u_{\varepsilon}\to u$ in $L^1_{\rm loc}$ by \eqref{e3.24}. Hence it suffices to show that there is $C$ independent of ${\varepsilon}$ such that
\begin{equation}\label{e3.52aa}
\|u_{\varepsilon}\|=\int_{{\mathbb{R}}^d}|u_{\varepsilon}(x)|\Phi(x)dx\le C,\ \ \forall \,{\varepsilon}>0,\end{equation}where $\Phi\in C^2({\mathbb{R}}^d)$ is such that $1\le\Phi(x),\ \forall \,x\in{\mathbb{R}}^d,$ and
$$\Phi(x)\to+{\infty}\mbox{\ \ as }|x|\to{\infty},\ \nabla\Phi\in L^{\infty},\ \Delta\Phi\in L^{\infty}.$$(An example is $\Phi(x)\equiv(1+|x|^2)^\alpha$ with $\alpha\in\left(0,\frac12\right].$)

We fix such a function $\Phi$ and
   assume that
$$f\in L^1\cap L^{\infty},\ \ \|f\|=\int_{{\mathbb{R}}^d}\Phi(x)|f(x)|dx <{\infty}.$$
If we  multiply equation \eqref{e3.10} by ${\varphi}_\nu{\mathcal{X}}_\delta(\widetilde \beta_{\varepsilon}(u_{\varepsilon}))$, where    ${\varphi}_\nu(x)=\Phi(x)\exp(-\nu\Phi (x))$, $\nu>0$, and integrate over ${\mathbb{R}}^d$, we get, since ${\mathcal{X}}'_\delta\ge0,$
\begin{equation}\label{e3.50a}
\begin{array}{l}
\displaystyle \int_{{\mathbb{R}}^d}u_{\varepsilon}{\mathcal{X}}_\delta
(\widetilde \beta_{\varepsilon}(u_{\varepsilon})){\varphi}_\nu\,dx
\le-{{\lambda}} \displaystyle \int_{{\mathbb{R}}^d}\nabla\widetilde \beta_{\varepsilon}(u_{\varepsilon})
\cdot\nabla({\mathcal{X}}_\delta(\widetilde \beta_{\varepsilon}(u_{\varepsilon})){\varphi}_\nu)dx\\
\qquad+{{\lambda}} \displaystyle \int_{{\mathbb{R}}^d}D_{\varepsilon} b^*_{\varepsilon}(u_{\varepsilon})\cdot
\nabla({\mathcal{X}}_\delta(\widetilde \beta_{\varepsilon}(u_{\varepsilon})){\varphi}_\nu)dx
+\displaystyle \int_{{\mathbb{R}}^d}|f|{\varphi}_\nu dx
\\
\qquad\le-{{\lambda}}\displaystyle \int_{{\mathbb{R}}^d}
(\nabla\widetilde \beta_{\varepsilon}(u_{\varepsilon})
\cdot\nabla{\varphi}_\nu){\mathcal{X}}_\delta(\widetilde \beta_{\varepsilon}(u_{\varepsilon}))dx\\
\qquad
+{{\lambda}} \displaystyle \int_{{\mathbb{R}}^d}D_{\varepsilon} b^*_{\varepsilon}(u_{\varepsilon})\cdot\nabla
\widetilde \beta_{\varepsilon}(u_{\varepsilon}){\mathcal{X}}'_\delta(\widetilde \beta_{\varepsilon}(u_{\varepsilon})){\varphi}_\nu dx\\
\qquad+{{\lambda}} \displaystyle \int_{{\mathbb{R}}^d}(D_{\varepsilon} \cdot\nabla{\varphi}_\nu)b^*_{\varepsilon}(u_{\varepsilon}){\mathcal{X}}_\delta(\widetilde \beta_{\varepsilon}(u_{\varepsilon}))dx
+\displaystyle \int_{{\mathbb{R}}^d}|f|{\varphi}_\nu dx.\end{array}\end{equation}
Here, $b^*_{\varepsilon}(u)=b_{\varepsilon}(u)u.$
Letting $\delta\to0$, we get as above
\begin{equation}\label{e3.50b}
\begin{array}{l}
\displaystyle \int_{{\mathbb{R}}^d}|u_{\varepsilon}|{\varphi}_\nu dx
\le-{{\lambda}}\displaystyle \int_{{\mathbb{R}}^d}\nabla|\widetilde \beta_{\varepsilon}(u_{\varepsilon})|\cdot\nabla{\varphi}_\nu dx\vsp
\quad
+\overline{\lim\limits_{\delta\to0}}\displaystyle \frac{{\lambda}}\delta
\displaystyle \int_{[|\widetilde \beta_{\varepsilon}(u_{\varepsilon})|\le\delta]}
|D_{\varepsilon}|\,|b^*_{\varepsilon}(u_{\varepsilon})|\,|\nabla\widetilde \beta_{\varepsilon}(u_{\varepsilon})|{\varphi}_\nu dx\vsp
\quad
+{{\lambda}}\displaystyle \int_{{\mathbb{R}}^d}({\rm sign}\,u_{\varepsilon}) b^*_{\varepsilon}(u_{\varepsilon})
(D_{\varepsilon}\cdot\nabla{\varphi}_\nu)dx+
\displaystyle \int_{{\mathbb{R}}^d}\!|f|{\varphi}_\nu dx\vsp\quad
\le{{\lambda}}\displaystyle \int_{{\mathbb{R}}^d}
(|\widetilde \beta_{\varepsilon}(u_{\varepsilon})|
\Delta{\varphi}_\nu
+|b^*_{\varepsilon}(u_{\varepsilon})|
\,|D_{\varepsilon}\cdot\nabla{\varphi}_\nu|)dx
+ \displaystyle \int_{{\mathbb{R}}^d} |f|{\varphi}_\nu dx,
\end{array}\end{equation}where in the last step we used that
$$|b^*_{\varepsilon}(u_{\varepsilon})|\le{\rm Lip}(b^*_{\varepsilon})|u_{\varepsilon}|\le\frac1{\varepsilon}\,{\rm Lip}(b^*_{\varepsilon})|\widetilde \beta_{\varepsilon}(u_{\varepsilon})|.$$
  We have
\begin{equation*}
\begin{array}{l}
\nabla{\varphi}_\nu(x)=(\nabla \Phi-\nu\Phi\nabla\Phi)\exp(-\nu\Phi),\vsp
\Delta{\varphi}_\nu(x)=(\Delta\Phi-\nu|\nabla\Phi|^2
-\nu\Phi\Delta\Phi+\nu^2\Phi|\nabla\Phi|^2
-\nu|\nabla\Phi|^2)\exp(-\nu\Phi).
\end{array}\end{equation*}
Then, letting $\nu\to0$ in \eqref{e3.50b}, since $M:=\sup\limits_{{\varepsilon}>0}|u_{\varepsilon}|_{\infty}<{\infty},$ $|b^*_{\varepsilon}(r)|\le|b|_{\infty}|r|$, $|D_\varepsilon|\le|D|$ and $|\widetilde \beta_{\varepsilon}(r)|\le\Big(\sup\limits_{|r|\le M}\beta'(r)+{\varepsilon}\Big)|r|,$ $\forall  r\in[-M,M]$, we get
$$  	\|u_{\varepsilon}\|
\le\|f\|+ C {{\lambda}}(|\Delta\Phi|_{\infty}+|D|_{\infty}|\nabla\Phi|_{\infty}) |f|_1,\ \forall {\varepsilon}\in(0,1).$$
Hence \eqref{e3.52aa} and, therefore, \eqref{e3.52a}, hold for all $f\in L^1\cap L^{\infty}$ with $\|f\|<{\infty}.$ Since the latter set is dense in $L^1$, we get \eqref{e3.52a}, as claimed.

\bigskip \noindent{\bf Proof of Theorem \ref{t2.2} (continued).} As seen earlier, the solution $u_h$ to the finite difference scheme \eqref{e2.2}--\eqref{e2.4} is uniformly convergent on every compact interval $[0,T]$ to $u\in C([0,{\infty});L^1)$. By \eqref{e3.3b} and \eqref{e3.8b}, by a standard argument we obtain that, for $u_0\in L^1\cap L^{\infty}$,
$$|u(t)|_{\infty}=|S(t)u_0|_{\infty}\le
\exp(|({\rm div}\, D)^-+|D||^{\frac12}_{\infty}t)|u_0|_{\infty},\ \forall  t\ge0.$$
\eqref{e2.5} follows by \eqref{e3.8b} and \eqref{e3.3bb}.

Let us prove now that $u$ is a distributional solution to \eqref{e1.1}. We note first that by \eqref{e2.4} we have (setting $u_h(t)=u_0$ for $t\in(-{\infty},0)$)
\begin{equation}\label{e3.33a}
\begin{array}{l}
u_h(t)-h\Delta \beta(u_h(t))+h\ {\rm div}(D b(u_h(t))u_h(t))=u_h(t-h),\ t\ge0,\vsp
u_h(0)=u_0.\end{array}
\end{equation}
Since $\lim\limits_{h\to0}u_h(t)=S(t)u_0$ in $L^1$ locally uniformly in $t\in[0,{\infty})$, we have by \eqref{e2.4a} that
$$|u_h(t)|_{\infty}\le\exp(|({\rm div}\, D)^-+|D||^{\frac12}_{\infty}t)|u_0|_{\infty}+1,\ t\ge0,$$for small enough $h$, and, hence, for $h\to0$, $\beta(u_h(t))\to\beta(u(t))$ in $L^1_{\rm loc},$ a.e. $t>0.$

Let ${\varphi}\in C^{\infty}_0([0,{\infty})\times{\mathbb{R}}^d)$.   Then,
  by \eqref{e3.33a} we have
$$\begin{array}{r}
\displaystyle \int^{\infty}_0\int_{{\mathbb{R}}^d}\frac 1h\ (u_h(t,x)-u_h(t-h,x)){\varphi}(t,x)
-\beta(u_h(t,x))\Delta{\varphi}(t,x)\\
-b(u_h(t,x))u_h(t,x)D_{\varepsilon}(x)\cdot\nabla{\varphi}(t,x)dt\,dx
=0,\end{array}$$if we take $u_h(t,x)\equiv u_0(x)$ for $t\in(-h,0]$.
Then, replacing the first term by
$$\displaystyle \int^{\infty}_0\int_{{\mathbb{R}}^N }\frac1h\,u_h(t,x)({\varphi}(t+h,x)-{\varphi}(t,x))dt\,dx+\frac1h\displaystyle \int^h_0\int_{{\mathbb{R}}^N}u_0(x){\varphi}(t,x)dt\,dx$$
and, letting $h\to0$, we get \eqref{e2.7}, as claimed. This completes the proof.\hfill$\Box$

\begin{remark}\label{r3.3}\rm Theorem \ref{t2.2} extends by a slight modification of the proof to the  multivalued functions $\beta:{\mathbb{R}}\to2^{\mathbb{R}}$ with $D(\beta)=(-\infty,+\infty)$ and which are maximal monotone graphs on ${\mathbb{R}}\times{\mathbb{R}}$, that is,
	 $(v_1-v_2)(u_1-u_2)\ge0,$ $ \forall  v_i\in\beta(u_i),$ $ i=1,2,$ and $R(1+\beta)={\mathbb{R}}.$
\end{remark}

 \section{Regularizing effect on initial data}\label{s5}
\setcounter{equation}{0}

Consider here equation \eqref{e1.1} under the following hypotheses.
\begin{itemize}
	\item[(k)] $\beta\in C^2({\mathbb{R}}),\ \beta'(r)\ge a|r|^{\alpha-1},\ \forall  r\in{\mathbb{R}};\
	\beta(0)=0,$\\  where $\alpha\ge1,$    $ d\ge3,$ $a>0.$
	\item[(kk)] $D\in L^{\infty}({\mathbb{R}}^d;{\mathbb{R}}^d),$ ${\rm div}\ D\in L^2(\mathbb{R}^d)+L^{\infty}({\mathbb{R}}^d)$, ${\rm div}\ D\ge0,$ a.e.
	\item[(kkk)] $b\in C_b({\mathbb{R}})\cap C^1({\mathbb{R}}),\ b\ge0.$
\end{itemize}
 \newpage

\begin{theorem}\label{t5.1} Let $d\ge3.$ Then, under Hypotheses {\rm(k), (kk), (kkk)}, the generalized solution $u$ to \eqref{e1.1} given by Theorem {\rm\ref{t2.2}} for $\mu=u_0 dx,$ $u_0\in L^1$, satisfies
\begin{equation}\label{e5.1}
u(t)\in L^{\infty},\ \forall  t>0, \end{equation}
\begin{equation}\label{e5.2}
|u(t)|_{\infty}\le C\,t^{-\frac d{2+(\alpha-1)d}}\ |u_0|^{\frac2{2+d(\alpha-1)}}_1,\ \forall  t\in(0,{\infty}),\ u_0\in L^1,
\end{equation}
where $C$ is independent of $u_0$.
\end{theorem}

\noindent{\bf Proof.} We shall first prove the following lemma.
\begin{lemma}\label{l5.2} Let $u_{{\lambda}}=(I+{{\lambda}} A)^{-1}f$, where $A$ is the operator \eqref{e3.5a}, \eqref{e3.5aa}. Then, for each $p>1$ and ${{\lambda}}>0$, we have	
\begin{equation}\label{e5.3}
|u_{{\lambda}}|^p_p+{{\lambda}} a C\frac{p(p-1)}{(p+\alpha-1)^2}\, |u_{{\lambda}}|^{p+\alpha-1}_{\frac{(p+\alpha-1)d}{d-2}}\le|f|^p_p,\ \forall  f\in L^p\cap L^1.\end{equation}
where $C$ is independent of $p$ and ${{\lambda}}$.
\end{lemma}

\noindent{\bf Proof.} Let $p\in(1,{\infty})$. By approximation, we may assume that $f\in L^1\cap L^{\infty}$. Let us first explain the proof by heuristic computations in the case $b\equiv1,$  that is, for the equation
  \begin{equation}\label{e5.4}
  u_{{\lambda}}-{{\lambda}}\Delta\beta(u_{{\lambda}})+{{\lambda}}\ {\rm div}(Du_{{\lambda}})=f .\end{equation}We multiply \eqref{e3.10} by $|u_{{\lambda}}|^{p-2}u_{{\lambda}}$ and integrate over  ${\mathbb{R}}^d$ and   get
 \begin{equation}\label{e4.4prim}\begin{array}{l}
  |u_{{\lambda}}|^p_p+(p-1){{\lambda}}\displaystyle \int_{{\mathbb{R}}^d}\beta'(u_{{\lambda}})|\nabla u_{{\lambda}}|^2|u_{{\lambda}}|^{p-2}dx\vsp
   \qquad\qquad=\displaystyle \int_{{\mathbb{R}}^d}f|u_{{\lambda}}|^{p-2}u_{{\lambda}} dx-\displaystyle \frac{(p-1)}p\ {{\lambda}}\displaystyle \int_{{\mathbb{R}}^d}|u_{{\lambda}}|^p
  {\rm div}\ D\,dx\vsp
    \qquad\qquad\le|f|_p|u_{{\lambda}}|^{p-1}_p
  \le\displaystyle \frac1p\ |f|^p_p+\left(1-\frac1p\right)\,|u_{{\lambda}}|^p_p.\end{array}\end{equation}
  Taking into account (k)  yields
  $$|u_{{\lambda}}|^p_p+ap(p-1){{\lambda}}\int_{{\mathbb{R}}^d}|u_{{\lambda}}|^{p+\alpha-3}|\nabla u_{{\lambda}}|^2dx\le|f|^p_p.$$On the other hand,
  $$\begin{array}{ll}
  \displaystyle \int_{{\mathbb{R}}^d}|u_{{\lambda}}|^{p+\alpha-3}
  |\nabla u_{{\lambda}}|^2dx
  \!\!\!&\displaystyle =\left(\frac2{p+\alpha-1}\right)^2\int_{{\mathbb{R}}^d}
  \left|\nabla
  \left(|u_{{\lambda}}|^{\frac{p+\alpha-1}2}\right)\right|^2dx\vsp&\displaystyle
 \ge C\left(\frac2{p+\alpha-1}\right)^2\left(\displaystyle \int_{{\mathbb{R}}^d}
 |u_{{\lambda}}|^{\frac{(\alpha-1+p)d}{d-2}}\ dx\right)^{\frac{d-2}d}\end{array}$$
  	(by the Sobolev embedding theorem in ${\mathbb{R}}^d$). This yields
  	 $$|u_{{\lambda}}|^p_p+{{\lambda}} a C\,\frac{p(p-1)}{(p+\alpha-1)^2}\,|u_{{\lambda}}|^{\alpha-1+p}_{\frac{(\alpha-1+p)d}{d-2}}\le|f|^p_p,\ \forall {{\lambda}}>0,$$as claimed. 
  	
  	 Of course, the above argument is heuristic since e.g. $u_{{\lambda}}\not\in H^1$. To make the proof rigorous, we recall that by Lemma \ref{l3.2a} the solution $u:=u_{{\lambda}}$ to \eqref{e3.8} constructed in Lemma \ref{l3.2} is an $L^1$-limit of solutions $u_{\varepsilon},$ ${\varepsilon}>0$, to the approximating equations \eqref{e3.17} (with $A_{\varepsilon}$ as in \eqref{e3.25b}). So, we shall start with \eqref{e3.17} (instead of \eqref{e5.4}) and with its solution $u_{\varepsilon}$. Then we know by the proof of Lemma \ref{l3.1} (see \eqref{e3.11}) that $u_{\varepsilon},b(u_{\varepsilon})u_{\varepsilon}\in H^1\cap L^1\cap L^{\infty},$ $\widetilde \beta_{\varepsilon}(u_{\varepsilon})\in H^2.$  We~ have,  for all $r\in{\mathbb{R}}$,
  	
  	 \begin{equation}\label{e4.4'}\widetilde \beta'_{\varepsilon}(r)=\frac{\beta'(g_{\varepsilon}(r))}{(1+{\varepsilon}\beta')(g_{\varepsilon}(r))}+{\varepsilon}\ge h_{\varepsilon}(g_{\varepsilon}(r)),\end{equation}where
   	 \begin{equation}\label{e4.6'}
   	 h_{\varepsilon}(r)=\frac{a|r|^{\alpha-1}}{1+{\varepsilon} a|r|^{\alpha-1}},\  r\in{\mathbb{R}},\ \ \ g_{\varepsilon}=(I+{\varepsilon}\beta)^{-1}.\end{equation}   	
  	 Define ${\varphi}_\delta:{\mathbb{R}}\to{\mathbb{R}}$ by
  	 $${\varphi}_\delta(r)=(|r|+\delta)^{p-2}r,\ \ r\in{\mathbb{R}}.$$Then, ${\varphi}_\delta\in C^1_b$, $\lim\limits_{\delta\to0}{\varphi}'_\delta(r)=(p-1)|r|^{p-2}$ and ${\varphi}'_\delta(r)\!\ge\!\min(1,p{-}1)(|r|{+}\delta)^{p-2}$. Now, we multiply \eqref{e3.17} by ${\varphi}_\delta(u_{\varepsilon})$ and  obtain  	
  	 \begin{equation}
  	 \label{e4.4''}
  	 \begin{array}{l}
  	 \displaystyle \int_{{\mathbb{R}}^d}u_{\varepsilon}{\varphi}_\delta(u_{\varepsilon})dx+{{\lambda}}\int_{{\mathbb{R}}^d}\widetilde \beta'_{\varepsilon}(u_{\varepsilon})|\nabla u_{\varepsilon}|^2{\varphi}'_\delta(u_{\varepsilon})dx\\\qquad={{\lambda}}\displaystyle \int_{{\mathbb{R}}^d}(D_{\varepsilon}\cdot\nabla u_{\varepsilon}){\varphi}'_\delta(u_{\varepsilon})b_{\varepsilon}(u_{\varepsilon})u_{\varepsilon} dx+\displaystyle \int_{{\mathbb{R}}^d}f{\varphi}_\delta(u_{\varepsilon})dx.
  	 \end{array}\end{equation}
  	 Defining
  	 $$\psi(r)=\int^r_0{\varphi}'_\delta(s)b_{\varepsilon}(s)s\,ds,$$we see that $\psi\ge0$, hence the first integral in the right hand side of \eqref{e4.4''} is equal to
  	 \begin{equation}
  	 \label{e4.4'''}
  	 -\displaystyle \int_{{\mathbb{R}}^d}({\rm div}\,D_{\varepsilon})\psi(u_{\varepsilon})dx\le C_\varepsilon,\end{equation}
  where (see \eqref{e3.23secund})
  	 	$$C_\varepsilon=\frac12\,|b|_\infty|D|_\infty\sup_{\varepsilon\in(0,1)}|u_\varepsilon|_\infty\int_{\mathbb{R}^d}
   	{\bf1}_{|x|\ge\frac1\varepsilon}|u_\varepsilon|dx.$$
  	 Furthermore, we deduce from \eqref{e4.4'} that the second integral on the left hand side of \eqref{e4.4''} dominates
  	 \begin{equation}
  	 \label{e4.4iv}
  	 \begin{array}{ll}
  	 \displaystyle \int_{{\mathbb{R}}^d}h_{\varepsilon}(g_{\varepsilon}(u_{\varepsilon}))|\nabla u_{\varepsilon}|^2{\varphi}'_\delta(u_{\varepsilon})dx\!\!&
  	 = \displaystyle \int_{{\mathbb{R}}^d}|\nabla\psi_{{\varepsilon},\delta}(u_{\varepsilon})|^2dx \\
  	 &\ge
  	 C\left(\displaystyle \int_{{\mathbb{R}}^d}|\psi_{{\varepsilon},\delta}(u_{\varepsilon})|^{\frac{2d}{d-2}}dx\right)^{\frac{d-2}d},\end{array} \end{equation}where
  	 \begin{equation}\label{e4.9'}
  	 \psi_{{\varepsilon},\delta}(r)=\int^r_0\sqrt{h_{\varepsilon}(g_{\varepsilon}(s)){\varphi}'_\delta(s)}\ ds,\ \  r\in{\mathbb{R}}.\end{equation}(Here,   we used the Sobolev embedding.) Combining \eqref{e4.4''}--\eqref{e4.4iv} and letting $\delta\to0$, we obtain by
  	 Fatou's lemma and \eqref{e3.20a}
  	   	  \begin{equation}
  	 \label{e4.4v}
  	 |u_{\varepsilon}|^p_p+\frac{{{\lambda}} C}p
  	 \left(\displaystyle \int_{{\mathbb{R}}^d}|\psi_{\varepsilon}(u_{\varepsilon})|^{\frac{2d}{d-2}}dx\right)^{\frac{d-2}2}
  	 \le\frac1p\ |f|^p_p+\left(1-\frac1p\right)|u_{\varepsilon}|^p_p+C_\varepsilon,  \end{equation}where
  	\begin{equation}
  	\label{e4.9prim}\psi_{{\varepsilon}}(r)=\sqrt{p-1}
  	 \int^r_0\sqrt{h_{\varepsilon}(g_{\varepsilon}(s))|s|^{p-2}}  \ ds,\ \  r\in{\mathbb{R}}.\end{equation}Obviously, $\psi_{\varepsilon}(u_{\varepsilon}),\ {\varepsilon}>0,$ are equicontinuous, hence by \eqref{e3.20a} and \eqref{e3.40d}
  	 $$(\psi_{\varepsilon}(u_{\varepsilon}))^2\to a(p-1)\left(\frac2{p+\alpha-1}\right)^2|u|^{p+\alpha-1},\mbox{\ \ a.e. on }{\mathbb{R}}^d,$$and, since $u_\varepsilon\to u$ in $L^1$  by Lemma \ref{l3.2a}, 
  	  $\lim_{\varepsilon\to0}C_\varepsilon=0.$  Therefore, by Fatou's lemma, \eqref{e4.4v} implies \eqref{e5.3}.

  	\bigskip\noindent{\bf Proof of Theorem \ref{t5.1}.} We choose $p=p_n$, where $\{p_n\}$ are defined by
  	$$p_{n+1}=\frac d{d-2}\ (p_n+\alpha-1),\   p_0>1.$$
  	  	Then, by \eqref{e5.3}, we get
  	  	 $$|u_{{\lambda}}|^{p_n}_{p_n}+Ca
  	  	 \frac{p_n(p_n-1)}{(p_n+\alpha-1)^2}\,
  	  	 {{\lambda}}|u_{{\lambda}}|^{p_n+\alpha-1}_{p_{n+1}}\le|f|^{p_n}_{p_n},\ n=0,1,...$$ 
  	  	 Now, we apply Theorem  5.2  in A. Pazy \cite{14}, where ${\varphi}_n(u)=|u|^{p_n}_{p_n},$ $\beta_n=\frac{d-2}d$, $C_n=Ca\,\inf\limits_{p\in[p_0,{\infty})}\,
  	  	 \frac{p(p-1)}{(p+\alpha-1)^2},$ and conclude that (see Proposition 6.5 in \cite{14} or~\cite{15}) that
  	  	  	  \begin{equation}\label{5.5}
  	  	  	  |u(t)|_{\infty}\le C_{p_0}t^{-\frac d{2p_0+(\alpha-1)d}}\ |u_0|^{\frac{2p_0}{2p_0+d(\alpha-1)}}_{p_0},\ \forall  t>0,\ u_0\in L^{p_0}.\end{equation}
    	Define
  	\begin{equation}
  	\label{5.5prim}
  	C_{\alpha,d}:=\frac{d+2}{2d}+\sqrt{(\alpha-1)
  	\left(\alpha+\frac2d\right)+\left(\frac{d+2}{2d}\right)^2}.
  	\end{equation}
  	Note that, since $\alpha>1-\frac2d$, the value under the root is strictly bigger than $\left(\frac{d-2}{2d}\right)^2$, hence $C_{\alpha,d}>1.$

  	  	 \begin{lemma}\label{l5.3} Let $p_0\in(1,C_{\alpha,d})$. Then, for some constant $C_{p_0}>0$,
  	  	 	\begin{equation}
  	  	 	\label{5.6}
  	  	 	|u(t)|_{p_0}\le C_{p_0}\,t^{-\frac{p_0-\gamma}{p_0(\gamma+\alpha-1)}}
  	  	 	|u_0|^{\frac{\gamma(p_0+\alpha-1}{p_0(\gamma+\alpha-1)}}_1,\ \forall   t>0,\  u_0\in L^1\cap L^{p_0},
  	  	 	\end{equation}
  	  	 	\begin{equation}
  	  	 	\label{5.7}
  	  	 	\gamma=\frac{2p_0+(\alpha-1)d}{(p_0+\alpha-2)d+2}\in(0,1).
  	  	 	\end{equation}
  	  	 \end{lemma}
    	
    	 \noindent{\bf Proof.} We may assume that $u_0\in L^{\infty}.$ We shall use the approximating scheme \eqref{e3.5}-\eqref{e3.6}. By estimate \eqref{e5.3},  we have, for   ${{\lambda}}=h$,
    	 $$|u^{i+1}_h|^{p_0}_{p_0}
    	 +Ch|u^{i+1}_h|^{p_0+\alpha-1}_{\frac{(p_0+\alpha-1)d}{d-2}}
    	 \le|u^i_h|^{p_0}_{p_0},\ i=0,1,...$$
    	 By the summation by parts formula, this yields, for all $t>0$,
    	 	\begin{equation}
    	 \label{5.8}
    	 t|u_h(t)|^{p_0}_{p_0}+ C\int^t_0
    	 s|u_h(s)|^{p_0+\alpha-1}_{\frac{(p_0+\alpha-1)d}{d-2}}\ ds
    	 	\le\int^t_0|u_h(s)|^{p_0}_{p_0}ds+h|u_0|^{p_0}_{p_0},
    	 	\end{equation}where $u_h$ is given by \eqref{e3.6} and where, here and below, by $C$ we denote various  constants independent of $t$ and $u_0$, but depending on $p_0$.

On the other hand,   by H\"older's inequality we have
$$|u_h(s)|^{p_0}_{p_0}\le|u_h(s)|_1^\gamma|u_h(s)|^{p_0-\gamma}_{\frac{(p_0+\alpha-1)d}{d-2}},\  s>0.$$
Then, substituting into \eqref{5.8}, we get 
   	 	    	 	
$$\begin{array}{l}
    	 	t|u_h(t)|^{p_0}_{p_0}+C\displaystyle \int^t_0s|u_h(s)|^{p_0+\alpha-1}_{\frac{(p_0+\alpha-1)d}{d-2}}\ ds-h|u_0|^{p_0}_{p_0}\vsp
    	 	\qquad\le|u_0|^\gamma_1\displaystyle \int^t_0|u_h(s)|^{p_0-\gamma}_{\frac{(p_0+\alpha-1)d}{d-2}}\ ds\vsp
    	 	\qquad\le|u_0|^\gamma_1\displaystyle \int^t_0
    	 	s^{\frac{p_0-\gamma}{p_0+\alpha-1}}
    	 	|u_h(s)|_{\frac{(p_0+\alpha-1)d}{d-2}})^{p_0-\gamma}\,s^{\frac{\gamma-p_0}{p_0+\alpha-1}}ds\vsp
    	 	\qquad\le|u_0|^\gamma_1
    	 	\left(\displaystyle \int^t_0
    	 	s|u_h(s)|^{p_0+\alpha-1}_{\frac{(p_0+\alpha-1)d}{d-2}}\ ds\right)^{\frac{p_0-\gamma}{p_0+\alpha-1}}
    	 	\left(\displaystyle \int^t_0
    	 	s^{\frac{\gamma-p_0}{\gamma+\alpha-1}}\,ds\right)^{\frac{\gamma+\alpha-1}{p_0+\alpha-1}}. 	 		\end{array}$$
 Since $p_0<C_{\alpha,d}$, we know by Lemma A.1 in the Appendix  that
 $$\frac{\gamma-p_0}{\gamma+\alpha-1} >-1.$$
 Hence the above is dominated by
 $$C|u_0|^\gamma_1\,t^{\frac{2\gamma+\alpha-p_0-1}{p_0+\alpha-1}}
 	\left(\int^t_0s|u_h(s)|^{p_0+\alpha-1}_{\frac{(p_0+\alpha-1)d}{d-2}}\,ds\right)^{\frac{p_0-\gamma}{p_0+\alpha-1}}.$$  	 		
 This yields, for $t>0$,
 $$t|u_h(t)|^{p_0}_{p_0}+ C\int^t_0s|u_h(s)|^{p_0+\alpha-1}_{\frac{(p_0+\alpha-1)d}{d-2}}ds
 \le C|u_0|^{\frac{\gamma(p_0+\alpha-1)}{\gamma+\alpha-1}}_1\,t^{\frac{2\gamma+\alpha-p_0-1}{\gamma+\alpha-1}}+h|u_0|^{p_0}_{p_0}.$$
Hence, dropping the integral on the left hand side and letting $h\to0$, we obtain that
$$|u(t)|_{p_0}\le C_{p_0}|u_0|^{\frac{\gamma(p_0+\alpha-1)}{p_0(\gamma+\alpha-1)}}_1\
    	 	t^{\frac{\gamma-p_0}{p_0(\gamma+\alpha-1)}},\  t>0,$$and  \eqref{5.6} is proved.
    	 	
 \bigskip\noindent{\bf Proof of Theorem  \ref{t5.1} (continued).} By approximation, we may assume that $u_0\in L^1\cap L^{\infty}$. Combining estimates \eqref{5.5} and \eqref{5.6}, we get, for $t>0,$
    	 \begin{equation}\label{5.9}
    	 \hspace*{-4mm}	\begin{array}{ll}
    	 	|u(t)|_{\infty}\!\!\!
    	 	&\le C_{p_0}\left(\frac t2\right)^{-\frac d{2p_0+(\alpha-1)d}}\ \left|u\left(\frac t2\right)\right|^{\frac{2p_0}{2p_0+d(\alpha-1)}}_{p_0}\vsp
    	 	&\le C_{p_0}\ t^{-\left(\frac d{2p_0+(\alpha-1)d}+
    	 		\frac{(p_0-\gamma)2p_0}{p_0(\gamma+\alpha-1)(2p_0+d(\alpha-1))}\right)}\
   |u_0|^{\frac{\gamma(p_0+\alpha-1)2p_0}{p_0(\gamma+\alpha-1)(2p_0+d(\alpha-1))}}_1\vsp
   &=C_{p_0}\ t^{-\frac d{2+(\alpha-1)d}}\ |u_0|^{\frac 2{2+(\alpha-1)d}}_1,
 \end{array}\end{equation}where we used Lemma A.2 in the Appendix in the last step and where again the constant $C_{p_0}$ changes from line to line. Hence, Theorem \ref{t5.1} is proved.
 
 \begin{remark}\label{r4.4}\rm In the proof of Theorem \ref{t5.1}, we have used  the weaker condition $\alpha>1-\frac2d$. However, since $\beta\in C^1(\mathbb{R})$, we should assume in (k) that $\alpha\ge1$. We plan to remove the assumption $\beta\in C^1(\mathbb{R})$ in a future paper.\end{remark}

 \section{Equation \eqref{e1.1} with a measure as initial  datum}\label{s6}
  	  	\setcounter{equation}{0}  	  	
  	  	
  	  	Consider here equation \eqref{e1.1} with the initial data $u_0=\mu\in{\mathcal{M}}_b.$
  	  	
  	  	\begin{definition}\label{d6.1}\rm The function $u:[0,{\infty})\to{\mathcal{M}}_b$ is a distributional solution to \eqref{e1.1} if  	  		
  	  		\begin{equation}
  	  		\label{e6.1}
  	  		u,\beta(u)\in L^1_{\rm loc}([0,{\infty})\times{\mathbb{R}}^d),\end{equation}
\begin{eqnarray}
\label{e6.2}
&&\hspace*{-10mm}\displaystyle \int^{\infty}_0\int_{{\mathbb{R}}^d}
u(t,x)({\varphi}_t(t,x)+b(u(t,x))D(x)\cdot\nabla{\varphi}(t,x))\\
&&\hspace*{-10mm}\quad+\beta(u(t,x))\Delta{\varphi}(t,x)dt\,dx+\mu({\varphi}(0,\cdot))=0,
 \ \forall {\varphi}\in C^{\infty}_0([0,{\infty})\times{\mathbb{R}}^d).\nonumber\end{eqnarray}
\end{definition}
We have
\begin{theorem}\label{t6.1}
Assume that Hypotheses {\rm(k), (kk), (kkk)} from Section {\rm\ref{s5}} hold and, in addition, that $D\in L^2(\mathbb{R}^d;\mathbb{R}^d)$ and
\begin{equation}\label{e5.2a}
|\beta(r)|\le C|r|^\alpha,\ \forall r\in\mathbb{R}.\end{equation}
Let $\mu\in{\mathcal{M}}_b.$
  Then, \eqref{e1.1} has a distributional solution which satisfies, for  \mbox{$dt$-a.e.} $t\in(0,{\infty})$,
  \begin{eqnarray}
  \label{e6.3prim}
 & u(t,x)\ge0,\mbox{ a.e. on }{\mathbb{R}}^d,\mbox{ provided }\mu\ge0,\\
 &\displaystyle
  \label{e6.3secund}
  \int_{{\mathbb{R}}^d}u(t,x)dx=\int_{{\mathbb{R}}^d}d\mu,\\
&\label{e6.4}
|u(t)|_{\infty}\le C\,t^{-\frac d{2+(\alpha-1)d}}\ \|\mu\|^{\frac2{2+d(\alpha-1)}}_{{\mathcal{M}}_b},
 \\
 &\label{e6.4aa}
 |u(t)|_1\le \|\mu\|_{{\mathcal{M}}_b}.
 \end{eqnarray}
 Furthermore, for all $p\in\left[1,\alpha+\frac2d\right),$
 \begin{eqnarray}
 \label{e6.5prim}
 u&\in& L^p((0,T)\times{\mathbb{R}}^d),\ \forall  T>0,\\
 \label{e6.5secund}
 \beta(u)&\in& L^1((0,T)\times{\mathbb{R}}^d),\ \forall  T>0.\end{eqnarray}
The map $t\mapsto u(t,x)dx\in{\mathcal{M}}_b$ has a $\sigma({\mathcal{M}}_b,C_b)$-continuous version on $(0,{\infty})$, denoted by $S(t)\mu,$ $t>0,$ for which  \eqref{e6.3prim}, \eqref{e6.3secund}, \eqref{e6.4} and \eqref{e6.4aa} hold for all $t>0$. Furthermore,
\begin{equation}
\label{e6.4a}
\lim_{t\to0}\int_{{\mathbb{R}}^d} (S(t)\mu)(x)\psi(x)dx=\mu(\psi),\ \forall \psi\in C_b.
\end{equation}
Defining $S(0)\mu=\mu$, then $S(t),$ $t\ge0,$  restricted to $L^1$ coincides with the semigroup from Theorem~{\rm\ref{t2.2}} and we have
$$
\|S(t)\mu_1-S(t)\mu_2\|_{{\mathcal{M}}_b}\le\|\mu_1-\mu_2\|_{{\mathcal{M}}_b},\ \ \forall \, t\ge0,\ \mu_1,\mu_2\in{\mathcal{M}}_b.$$\end{theorem}

\noindent{\bf Proof.}   Consider a smooth approximation $\mu_{\varepsilon}$ of $ u_0=\mu$ of the form
$$\mu_{\varepsilon}(x)=(\mu*\rho_{\varepsilon}),\ {\varepsilon}>0,$$where $\rho_{\varepsilon}(x)=\frac1{\varepsilon}\ \rho\left(\frac1{\varepsilon}\ |x|\right),\ \rho\in C^{\infty}_0([-1,1]),$ $\int^1_{-1}\rho(r)\, dr=1.$ Then, by Theorem \ref{t5.1}, the equation
\begin{equation}
\label{e6.5}
\begin{array}{l}
(u_{\varepsilon})_t-\Delta\beta(u_{\varepsilon})+{\rm div}(Db(u_\varepsilon)u_{\varepsilon})=0\mbox{ in }(0,{\infty})\times{\mathbb{R}}^d,\vsp
u_{\varepsilon}(0)=\mu_{\varepsilon},\end{array}\end{equation}has, for each ${\varepsilon}>0$, a unique generalized solution $u_{\varepsilon}\in C([0,{\infty});L^1)\cap L^{\infty}((\delta,{\infty})\times{\mathbb{R}}^d)$, $\forall \delta>0.$ More precisely, we have 
 \begin{equation}
 \label{e6.6}
  |u_{\varepsilon}(t)|_{\infty}\le C\, t^{-\frac d{2+(\alpha-1)d}}\ |\mu_{\varepsilon}|_1^{\frac2{2+d(\alpha-1)}}
 \le C\,t^{-\frac d{2+(\alpha-1)d}}\ \|\mu\|_{{\mathcal{M}}_b}^{\frac2{2+d(\alpha-1)}},\ t>0.\end{equation}
   Everywhere in the following, $C$ is a positive constant independent of $t$ and $\mu$ possibly changing from line to line. Also, for simplicity, we set $\|\mu\|=\|\mu\|_{{\mathcal{M}}_b}$, and     $\|\mu\|^{\frac2{2+(\alpha-1)d}}\,t^{-\frac d{2+(\alpha-1)d}}=\nu(t,\mu),\ \forall  t>0,\ \mu\in{\mathcal{M}}_b.$
 We also have by \eqref{e2.8}
 \begin{equation}
 \label{e6.7}
 |u_{\varepsilon}(t)|_1\le|\mu_{\varepsilon}|_1\le\|\mu\|,\ \forall  t\ge0,\ \forall {\varepsilon}>0.\end{equation}
 If we  formally multiply \eqref{e6.5} by $\beta(u_{\varepsilon})$ and integrate over $(\delta,t)\times{\mathbb{R}}^d$, for $\psi(r)\equiv\int^r_0\beta'(s)b(s)s\,ds,$   since $\psi,{\rm div}\,D\ge0$, we get
  \begin{equation}
  \label{e6.8}
  \begin{array}{l}
  \displaystyle \int_{{\mathbb{R}}^d}g(u_{\varepsilon}(t,x))dx+\int^t_\delta\int_{\mathbb{R}^d}|\nabla\beta (u_{\varepsilon}(s))|^2_2dx\,ds\\
  \qquad\quad=\displaystyle \int^t_\delta\int_{{\mathbb{R}}^d}\nabla(\psi(u_{\varepsilon}))\cdot D\,dx\,ds +\int_{{\mathbb{R}}^d}g(u_{\varepsilon}(\delta,x))dx\vsp
  \qquad\quad\le \displaystyle \int_{{\mathbb{R}}^d}g(u_{\varepsilon}(\delta,x))dx\le  C\|\mu\|(\nu(\delta,\mu))^{\alpha},\ \forall  t>\delta,\end{array}
  \end{equation}
  where  $g(r)\equiv\int^r_0\beta(s)ds\ge0$ and the last inequality follows by \eqref{e5.2a}. 
    Estimate \eqref{e6.8} can be derived rigorously by using   the finite dif\-fe\-rence  scheme \eqref{e3.5}--\eqref{e3.6} corresponding to the resolvent of the regularized version \eqref{e3.17} of  equation  \eqref{e6.5}. Indeed, by Lemma \ref{l3.2a}, it follows via the Trotter--Kato theorem for nonlinear semigroups (see, e.g., \cite{1}, p.~168) that, for each ${\varepsilon}>0$,
  $$u_{\varepsilon}(t)=\lim_{\nu\to0}\lim_{n\to{\infty}}\left(I+\frac tn\,A_\nu\right)^{-n}\mu_{\varepsilon},$$where $A_\nu$ is the operator defined by \eqref{e3.25b} and both limits are in $L^1$, locally uniformly in $t\in[0,{\infty})$. Hence,
  \begin{equation}
  \label{e6.14a}u_{\varepsilon}(t)=\lim_{\nu\to0}\lim_{h\to0} u_{\nu,h}(t),\ t\in[0,T],
  \end{equation}
  where
 \begin{equation}
\label{6.15a}
\begin{array}{l}
  u_{\nu,h}(t)=u^{i+1}_{\nu,h},\ t\in(ih,(i+1)h],\vsp
u^{i+1}_{\nu,h}+hA_\nu u^{i+1}_{\nu,h}=u^i_{\nu,h},\ i=0,1,...,N-1;\ \ Nh=T,\ u^0_{\nu,h}=\mu_{\varepsilon}.\end{array}\end{equation}

We know by the proof of Lemma \ref{l3.1} that, if $v\in L^1\cap L^{\infty}$, then for the solution $u_h$ to the equation $u_h+hA_\nu u_h=v$  (see \eqref{e3.25b} and \eqref{e3.17}), we have
$u_h,b_\nu(u_h)u_h\in H^1\cap L^1\cap L^{\infty},$ $\widetilde \beta_\nu(u_h)\in H^2$ and $|u_h|_{\infty}\le|v|_{\infty}.$   Hence, if we  multiply \eqref{6.15a} by $\widetilde \beta_\nu(u^{i+1}_{\nu,h})$ and integrate over ${\mathbb{R}}^d$, we get as above 
$$\begin{array}{l}
\displaystyle \int_{{\mathbb{R}}^d}g_\nu(u^{i+1}_{\nu,h}(x))dx+h\int_{{\mathbb{R}}^d}|\nabla\widetilde \beta_\nu(u^{i+1}_{\nu,h})|^2dx\\
\qquad\le
\displaystyle \int_{{\mathbb{R}}^d}g_\nu(u^i_{\nu,h})dx,\ i=0,1,...,N_1,\ Nh=T,\end{array}$$
where $g_\nu(r)=\int^r_0\widetilde \beta_\nu(s)ds$ and we used that ${\rm div}\,D_\nu\ge0,$ since $D_\nu=D$ because $D\in L^2(\mathbb{R}^d;\mathbb{R}^d).$ Summing over from $j=[N\delta/T]+1$ to $k-1=[Nt/T]$, we~get
$$\begin{array}{l}
\displaystyle \int_{{\mathbb{R}}^d}g_\nu(u^{k}_{\nu,h})dx+\frac h2\sum^{k-1}_{i=j}\int_{{\mathbb{R}}^d}|\nabla
\widetilde \beta_\nu(u^{i+1}_{\nu,h})|^2dx \le
 \displaystyle \int_{{\mathbb{R}}^d}g_\nu(u^j_{\nu,h})dx,\ \forall  k.\end{array}$$
Then, letting $h\to0$ and afterwards $\nu\to0$,
by \eqref{e6.14a} and, since $|u^i_{\nu,h}|_{\infty}\le|\mu_{\varepsilon}|_{\infty}$,  the closedness of the gradient on $L^2(0,T;L^2)$ and the weak lower semicontinuity implies   \eqref{e6.8}, as claimed.

Multiplying \eqref{e6.5} by $|u_{\varepsilon}|^{q-2}u_{\varepsilon}$, $q\ge2,$ and integrating over $(\delta,t)\times{\mathbb{R}}^d$, we~get by hypothesis (k) that
 \begin{equation}\label{e6.16*}
 \begin{array}{r}
  a(q-1)\displaystyle\left(\frac2{q{+}\alpha{-}1}\right)^2\!\!\! \int^t_\delta\!\!\int_{{\mathbb{R}}^d}\!
  |\nabla|u_{\varepsilon}|^{\frac{q+\alpha-1}2}|^2 ds\,dx
  +\frac1q\int_{{\mathbb{R}}^d}|u_{\varepsilon}(t,x)|^qdx\vsp
 \quad
  \le\displaystyle \frac1q\int_{{\mathbb{R}}^d}|u_{\varepsilon}(\delta,x)|^qdx \le C\|\mu\|(\nu(\delta,\mu))^{q-1}.\end{array}\end{equation}
  As in the previous case, the above calculus can be made rigorous if we replace \eqref{e6.5} by its discrete version \eqref{6.15a}, which we multiply   by $|u^{i+1}_{\nu,h}|^{q-2} u^{i+1}_{\nu,h}$ and integrate over  ${\mathbb{R}}^d$.  Indeed, noting that, since $u^{i+1}_{\nu,h},b_\nu(u^{i+1}_{\nu,h})u^{i+1}_{\nu,h}\in H^1\cap L^1\cap L^{\infty}$ and  $\widetilde \beta'_\nu(u^{i+1}_{\nu,h})\in H^2$ (see \eqref{e3.11}, \eqref{e3.36prim}), we get similarly as in the proof of Lemma~\ref{l5.2}   for $\psi_\nu(r)=\sqrt{p-1}\int^r_0\sqrt{h_\nu(g_\nu(s))|s|^{q-2}}\,ds,\ r\in{\mathbb{R}},$ where $h_{\varepsilon},g_{\varepsilon}$ are as in~\eqref{e4.6'},
  $$\begin{array}{r}
 \displaystyle \frac1q\int_{{\mathbb{R}}^d}|u^{i+1}_{\nu,h}|^qdx+h \int_{{\mathbb{R}}^d}
 |\nabla \psi_\nu(u^{i+1}_{\nu,h})|^2  dx \le\displaystyle \frac1q\int_{{\mathbb{R}}^d}|u^{i}_{\nu,h}|^qdx,\vsp i=0,...,N-1,\ Nh=T.\end{array}$$
 Summing over $i$ from $j=[N\delta/T]+1$ to $k-1=[Nt/T]$, we get

 $$\begin{array}{r}
  \displaystyle \frac1q\int_{{\mathbb{R}}^d}|u^{k}_{\nu,h}|^qdx+h
  \sum^{k-1}_{i=j}\int_{{\mathbb{R}}^d}
  |\nabla \psi_\nu(u^{i+1}_{\nu,h})|^2 dx
   \le\displaystyle \frac1q\int_{{\mathbb{R}}^d}|u^{j}_{\nu,h}|^qdx,\\ i=0,...,N-1,\ Nh=T.\end{array}$$
Letting $ h\to0$, and afterwards $\nu\to0$, \eqref{e6.16*} follows from \eqref{e6.14a} and the closedness of the gradient on $L^2(0,T;L^2).$

   Now, taking into account that by \eqref{e6.16*}, with $q=2p-1-\alpha$,  we get
    \begin{equation}\label{e6.16}
       \displaystyle
    \int^t_\delta|\nabla(|u_{\varepsilon}|^{p-1})|^2_2ds\le \|\mu\|(\nu(\delta,\mu))^{2(p-1)-\alpha}),\ \ \forall  t\ge\delta, \forall  p\ge\frac{\alpha+3}2.\end{equation}
    Moreover, by \eqref{e6.8},   $\{\nabla\beta(u_{\varepsilon})\}_{{\varepsilon}>0}$ is bounded in $L^2(\delta,T,L^2)$, and so
   \begin{equation}\label{e6.16a}
  \|\Delta\beta(u_{\varepsilon})-{\rm div}(Db(u_{\varepsilon})u_{\varepsilon})\|_{L^2(\delta,T;H^{-1})}\le C,\ \forall {\varepsilon}>0.\end{equation}
   Note also that, by \eqref{e6.16}, it follows that
  \begin{equation}\label{e5.19'}
  \int^t_\delta|\nabla(|u_{\varepsilon}|^{p-2}u_{\varepsilon})|^2_2ds\le C,\ \forall {\varepsilon}>0.\end{equation}
   Hence,  $\{|u_{\varepsilon}|^{p-1}\}_{{\varepsilon}>0}$ is bounded in $L^2(\delta,T;H^1)$ and so, by \eqref{e6.16a}, we infer that, for $m\ge4$,
   $$\|\,|u_{\varepsilon}|^{p-1}(\Delta\beta(u_{\varepsilon})\|_{L^1(\delta,T;H^{-m})}+\|{\rm div}(Db(u_{\varepsilon})u_{\varepsilon}))\|_{L^2(\delta,T;H^{-1})}\le C.$$This implies that the set
   $$\left\{\frac{\partial}{{\partial} t}(|u_{\varepsilon}|^{p-1}u_{\varepsilon})\right\}_{{\varepsilon}>0}
=\left\{p|u_{\varepsilon}|^{p-1}(\Delta\beta(u_{\varepsilon})-{\rm div}(Db(u_{\varepsilon})u_\varepsilon))\right\}_{{\varepsilon}>0}$$is bounded in $L^1(\delta,T;H^{-1}).$ Note that by \eqref{e6.16} applied to $p+1$ replacing $p$, we have that also $\{|u_{\varepsilon}|^{p-1}u_{\varepsilon}\}_{{\varepsilon}>0}$ is bounded in $L^2(\delta,T;H^{1})$.

  Then, by the Aubin-Lions-Simon compactness theorem (see \cite{17a}), the set $\{|u_{\varepsilon}|^{p-1}u_{\varepsilon}\}_{{\varepsilon}>0}$   is relatively compact in $L^2(\delta,T;L^2_{\rm loc})$ for all $0<\delta<T<{\infty}.$ Hence, along a subsequence, we have for $\gamma(r):=|r|^{p-1}r,$ $r\in{\mathbb{R}}$,
  	\begin{equation}
  \label{e6.11}
  \gamma(u_{\varepsilon})\to v,\ \mbox{ a.e. on }(0,{\infty})\times{\mathbb{R}}^d.\end{equation}
  Then, since $\gamma$ has a continuous inverse and since $\beta$ is continuous,  we have

\begin{equation}
\label{e5.20'}
u_{\varepsilon}\to u=\gamma^{-1}(v)\mbox{ and }  \beta(u_{\varepsilon})\to\beta(u),\mbox{ a.e. on }(0,{\infty})\times{\mathbb{R}}^d.\end{equation}
By   \eqref{e6.6}, \eqref{e6.7}, we have, for all $p\ge1,$   $\forall    t>0,$ $ {\varepsilon}>0,$
  \begin{equation}
  \label{e6.12}
  \begin{array}{r}
  |u_{\varepsilon}(t)|_{p}  \le|u_{\varepsilon}(t)|^{\frac1p}_{1}|u_{\varepsilon}(t)|_{{\infty}}^{\frac{p-1}p}
  \le\displaystyle  C\|\mu\|^{\frac{2p+(\alpha-1)d}{(2+(\alpha-1)d)p}}\,t^{-\frac{d(p-1)}{(2+(\alpha-1)d)p}}.  \end{array}
  \end{equation}
We have $\frac{d(p-1)}{2+(\alpha-1)d}<1$
  for every $p\in\left[1,\alpha+\frac2d\right)$. 
 Then, for such $p$,  \eqref{e6.12} implies that, for every $T>0$,
  \begin{equation}
  \label{e6.15a}
  \int^T_0|u_{\varepsilon}(t)|^p_{p}dt\le C,\end{equation}
and, therefore,  if  in addition  $p>1$, along a subsequence ${\varepsilon}\to0$,
   \begin{equation}
  \label{e6.14}
  u_{\varepsilon}\to u\mbox{ weakly in }L^p((0,T)\times {\mathbb{R}}^d).\end{equation}

  Moreover, by \eqref{e5.2a}, \eqref{e5.20'}  and \eqref{e6.15a}, it follows that $\{\beta(u_{\varepsilon})\}$ is bounded in $L^q((0,T)\times{\mathbb{R}}^d)$ for all $q\in\left(1,1+\frac2{\alpha d}\right)$, and so   (along a subsequence)
  \begin{equation}\label{e6.17a}
  \beta(u_{\varepsilon})\to\beta(u)\mbox{ weakly in $L^q((0,T)\times {\mathbb{R}}^d).$ }\end{equation}
Since, by \eqref{e2.7}, we have
  $$\int^{\infty}_0\int_{{\mathbb{R}}^d}(u_{\varepsilon}({\varphi}_t+D\cdot\nabla{\varphi})+\beta(u_{\varepsilon})\Delta{\varphi})dt\,dx
  +\int_{{\mathbb{R}}^d}\mu_{\varepsilon}(x){\varphi}(0,x)dx=0$$for any ${\varphi}\in C^{\infty}_0([0,{\infty})\times{\mathbb{R}}^d),$ letting ${\varepsilon}\to0$, we see by \eqref{e6.14}  and \eqref{e6.17a}    that
   $u$ satisfies \eqref{e6.2}. As regards
\eqref{e6.3prim},  \eqref{e6.3secund}, \eqref{e6.4} and \eqref{e6.4aa}, these by \eqref{e5.20'} immediately follow from the correspon\-ding properties of $u_{\varepsilon}$ and \eqref{e6.6}. Furthermore, \eqref{e6.5prim} follows from \eqref{e6.15a}
and Fatou's Lemma. Taking  $p=\alpha$ in \eqref{e6.5prim},   \eqref{e6.5secund} follows by \eqref{e5.2a}.

By \eqref{e6.3secund} and \eqref{e6.5secund}, we may apply Lemma 8.1.2 in \cite{0}, to conclude that $t\mapsto u(t,x)dx\in{\mathcal{M}}_b$ has a $\sigma({\mathcal{M}}_b,C_b)$-continuous version on $(0,{\infty})$,   denoted by $\mu_t,$ $t>0$.
To show \eqref{e6.4a}, we apply \eqref{e6.2} with ${\varphi}(t,x)=\psi(t)\zeta(x),$   $\psi\in C^{\infty}_0([0,{\infty}))$ and $\zeta\in C^{\infty}_0({\mathbb{R}}^d)$. Then, for
$$L\zeta(t,x)=\beta(u(t,x))\Delta\zeta(x)+D(x)\cdot\nabla\zeta(x),$$we have from \eqref{e6.2}
\begin{equation}
\label{e6.19}
\displaystyle \int^{\infty}_0\psi(t)\int_{{\mathbb{R}}^d}L\zeta\,d\mu_t\,dt
+\psi(0)\int_{{\mathbb{R}}^d}\zeta\,d\mu=-\displaystyle \int^{\infty}_0\frac d{dt}\,\psi(t)\int_{{\mathbb{R}}^d}\zeta\,d\mu_t\,dt,
\end{equation}
hence, choosing $\psi\in C^{\infty}_0((0,{\infty}))$, we obtain for $dt$-a.e. $t\in(0,{\infty})$,
\begin{equation}
\label{e6.20}
\int_{{\mathbb{R}}^d}\zeta\,d\mu_t =C+\int^t_0\int_{{\mathbb{R}}^d} L\zeta\,d\mu_s ds.
\end{equation}
By \eqref{e6.5secund}, the right hand side is continuous in $t\in[0,{\infty})$ and equal to $C$ at $t=0$, while,  as seen above, also the left hand side is continuous in $t\in(0,{\infty})$. Hence, we obtain that \eqref{e6.20} holds for all $t\in(0,{\infty})$ and
$$\lim_{t\to0} \int_{{\mathbb{R}}^d}\zeta\,d\mu_t=C.$$
Plugging \eqref{e6.20} into the right hand side of \eqref{e6.19}, with $\psi\in C^{\infty}_0([0,{\infty}))$ such that $\psi(0)=1$ and integrating by parts, we find
$$\int^{\infty}_0\psi\int_{{\mathbb{R}}^d}L\zeta\,d\mu_t\,dt+\int_{{\mathbb{R}}^d}\zeta\,d\mu=C+\int^{\infty}_0\psi\int_{{\mathbb{R}}^d}L\zeta\,d\mu_t\,dt$$and \eqref{e6.4a} follows, because  \eqref{e6.3secund} holds for all $t>0$, as we shall see below.
It is obvious that, for the $\sigma({\mathcal{M}}_b,C_b)$-continuous version $t\mapsto\mu_t$ of $t\mapsto u(t,x)dx$ on $(0,{\infty})$, properties \eqref{e6.3prim}, \eqref{e6.3secund} and \eqref{e6.4aa} hold for all $t>0$. For this version, it is also easily seen that
 $t\mapsto|u(t)|_{\infty}$ is lower semicontinuous, hence also \eqref{e6.4} follows for all $t>0$.

It remains to prove the last assertion in Theorem \ref{t6.1}. To express the dependence of our $\sigma({\mathcal{M}}_b,C_b)$-continuous version $[0,{\infty})\ni t\mapsto\mu_t\in{\mathcal{M}}_b$ with $\mu_0=\mu$ of our solution to \eqref{e6.2}, we set, for $\mu\in{\mathcal{M}}_b,$
 $P(t)\mu=\mu_t,\ \ t\ge0,$ and recall that $\mu_t$ has a density in $L^1$ for $t>0$, which we identify with $\mu_t$, i.e., $\mu_t\in L^1$, $\forall  t>0.$ Let $T>0$. By construction, we know that (along a subsequence depending on $\mu$) ${\varepsilon}\to0$
\begin{equation}
\label{e5.28}
S(\cdot)(\mu*\rho_{\varepsilon})\!\to P(\cdot)\!\mu,\mbox{ a.e. on }(0,T)\!\times\!{\mathbb{R}}^d\mbox{ and weakly in }L^p((0,T)\!\times\!{\mathbb{R}}^d)
\end{equation}
as functions of $(t,x)$  for $p\in\left(1,\alpha+\frac2d\right)$ (see \eqref{e5.20'}, \eqref{e6.14}, respectively). Here $S(t),$ $t\ge0$, is the semigroup from Theorem \ref{t6.1}.

\medskip\noindent{\bf Claim.} {\it If $\mu\in L^1$, then}
 $S(t)\mu=P(t)\mu \mbox{\ \ for all }t\ge0.$\medskip

\noindent To prove the claim  we recall that, since $\mu\in L^1$, we have  $\mu*\rho_{\varepsilon}\to\mu$ in $L^1$. Hence, by \eqref{e2.8} and \eqref{e5.28},
$$S(t)\mu=P(t)\mu\mbox{ in }{\mathcal{M}}_b\mbox{ for }dt-\mbox{a.e. }t\in[0,T].$$Since both sides are $\sigma(L^1,C_b)$-continuous in $t\in[0,T]$, this equality holds $\forall \,t\in[0,T]$, $T>0$, and the Claim   is proved.

Therefore, we may rename $P(t):{\mathcal{M}}_b\to{\mathcal{M}}_b$, $t\ge0$, and set $S(t)=P(t),$ $t\ge0$, since it is an extension of $S(t):L^1\to L^1$ for every $t\ge0$.

Finally, for $\mu,\widetilde \mu\in{\mathcal{M}}_b$ with corresponding solutions $u_{\varepsilon},\widetilde  u_{\varepsilon}$ to \eqref{e6.5}, we have by \eqref{e2.8}, for all $t\ge0$,
$$|u_{\varepsilon}(t)-\widetilde  u_{\varepsilon}(t)|_{1}\le|(\mu-\widetilde \mu)*\rho_{\varepsilon}|_{1}\le\|\mu-\widetilde \mu|_{{\mathcal{M}}_b}.$$Hence, for all ${\varphi}\in C_b([0,{\infty})),$ ${\varphi}\ge0,$ by \eqref{e5.20'} and Fatous's lemma, letting ${\varepsilon}\to0$  we get
$$\int^{\infty}_0{\varphi}(t)|S(t)\mu-S(t)\widetilde \mu|_{L^1}dt\le\int^{\infty}_0{\varphi}(t)\|\mu-\widetilde \mu\|_{{\mathcal{M}}_b}dt,$$so,
$$|S(t)\mu-S(t)\widetilde \mu|_{1}\le\|\mu-\widetilde \mu\|_{{\mathcal{M}}_b}\mbox{ for }dt-\mbox{a.e. }t\in(0,{\infty}).$$But the left hand side is lower semicontinuous in $t\in[0,{\infty})$, since $t\mapsto S(t)\mu$ is $\sigma(L^1,C_b)$-continuous, hence
$$\|S(t)\mu-S(t)\widetilde \mu\|_{{\mathcal{M}}_b}\le\|\mu-\widetilde \mu\|_{{\mathcal{M}}_b},\ \forall  \,t\in[0,{\infty}).$$

\begin{remark}
	\label{r6.3}\rm One might suspect that, if $\mu\ge0$, then under the hypotheses of Theorem \ref{t6.1} the solution $u$  is the unique nonnegative solution to \eqref{e6.2}. This is true for the porous media equation (\cite{13a}). If the uniqueness is true  for all $\mu\in L^1$, it follows by Theorem \ref{t2.2} that $u(t)\in C([\delta,{\infty});L^1)$, for each $\delta>0$. In~this case, our solutions in Theorem \ref{t6.1} starting from any $\mu\in {\mathcal{M}}_b$, would also have the flow property.
\end{remark}

\begin{remark}\label{r5.5}\rm It should be noted that, as follows from \eqref{e5.19'}, \eqref{e5.20'}, for  every $\mu\in{\mathcal{M}}_b$, $p\ge\frac{\alpha+3}2$,
	$$|S(\cdot)\mu|^{p-2}S(\cdot)\mu\in L^2(\delta,T;H^1),\ 0<\delta<T<{\infty}.$$
\end{remark}

\begin{remark}\label{r5.5a}\rm  If $D\equiv0$ and $\beta(r)\equiv a r^\alpha$, where $0<\alpha<\frac{d-2}2$, then (see \cite{19a}) equation \eqref{e1.1} has a nonnegative solution for $\mu\in{\mathcal{M}}_b({\mathbb{R}}^d)$ if and only if $\mu(K)=0$ for each compact set $K$ of $C_{2,p}$-capacity zero, where  $p=\frac1{1-\alpha}.$ The extension of this result to the present case remains to be done. (See also Remark \ref{r4.4}.)
	\end{remark}

\section{The McKean-Vlasov equation}\label{s4}
\setcounter{equation}{0}

As a direct consequence of Theorems \ref{t2.2} and \ref{t6.1}, we obtain (probabilistically) weak  solutions to the McKean-Vlasov SDE \eqref{e1.2}. More precisely, we have

\begin{theorem}\label{t7.1} Assume that we are in one of the following situations:  
	\begin{itemize}
		\item[\rm(a)] Hypotheses {\rm(i), (ii), (iii)} from Section \ref{s1} and \eqref{e2.4prim} hold and let $u$ be the solution of \eqref{e2.7} from Theorem {\rm\ref{t2.2}} with the initial condition $\mu=u_0dx$,  $u_0\in{\mathcal{P}}_0({\mathbb{R}}^d)\cap L^{\infty}$.
		\item[\rm(b)] Hypotheses {\rm(k), (kk), (kkk)} from Section {\rm\ref{s5}} and \eqref{e5.2a} hold and, addi\-tionally, let $D\in L^2(\mathbb{R}^d;\mathbb{R}^d)$.  Let $u$ be the solution of \eqref{e6.2} from Theorem {\rm\ref{t6.1}} with the initial condition   \mbox{$\mu\in{\mathcal{P}}({\mathbb{R}}^d)$.}
	\end{itemize} Then,
 there exists a $($probabilistically$)$ weak solution $X$ to \eqref{e1.2} on some filtered probability space $({\Omega},{\mathcal{F}},\mathbb{P},({\mathcal{F}}_t)_{t\ge0})$ with an ${\mathbb{R}}^d$-valued $({\mathcal{F}}_t)$-Brownian motion $W(t),$ $t\ge0$, such that $X$ has $\mathbb{P}$-a.s. continuous sample paths and
	\begin{equation}\label{e4.1}
	\mu=\mathbb{P}\circ(X(0))^{-1}\mbox{\ \ and\ \ }u(t,x)dx=\mathbb{P}\circ(X(t))^{-1}(dx),\ \forall  t>0.\end{equation}
	\end{theorem}

\noindent{\bf Proof.} By our assumptions, we have, for every $T>0$,
$$\int^T_0\int_{{\mathbb{R}}^d}(|\beta(u(t,x))|+|b(u(t,\cdot))|\,|D(x)|u(t,x))dx\,dt<{\infty}.$$Hence, the assertion follows immediately by Section 2 in \cite{5}.\hfill$\Box$

\section*{Appendix}

\noindent{\bf Lemma A.1.} {\it Let $\alpha\in\left(\frac{d-2}d,{\infty}\right)$, $p_0\in(1,{\infty})$. Let $C_{\alpha,d}$ and $\gamma$ be as defined in \eqref{5.5prim}, \eqref{5.7}, respectively. Then
	\begin{itemize}
		\item[\rm(i)] $\gamma=1-\displaystyle \frac{(p_0-1)(d-2)}{(p_0+\alpha-2)d+2},$   $\gamma+\alpha-1>0.$
	\item[\rm(ii)] If $p_0<C_{\alpha,d}$, then $\displaystyle \frac{\gamma-p_0}{\gamma+\alpha-1}>-1.$	
	\end{itemize}}

The proof follows by a direct computation and will be omitted.\smallskip

By Lemma A.1, it also follows

\medskip \noindent{\bf Lemma A.2.} {\it Consider the situation of Lemma {\rm A.1}. Then
\begin{itemize}
	\item[\rm(j)] $p_0-\gamma=\frac{(p_0-1)(p_0+\alpha-1)d}{(p_0+\alpha-2)d+2},$  $\gamma+\alpha-1=\frac{(2+(\alpha-1)d)(p_0+\alpha-1)}{(p_0+\alpha-2)d+2}.$
		\item[\rm(jj)] $\displaystyle \frac{2\gamma(p_0+\alpha-1)}{(\gamma+\alpha-1)(2p_0+(\alpha-1)d)}=\frac2{2+(\alpha-1)d}.$
\end{itemize}}

\bigskip \noindent{\bf Acknowledgement.} This work was supported by the DFG through CRC 1283. The authors are grateful to the anonymous referee for very useful  comments and suggestions.

\end{document}